\documentclass[11pt]{article}
\usepackage[dvips]{graphics}
\usepackage{multirow}
\usepackage{array}
\usepackage{epsfig,rotating}
\usepackage{amssymb,amsmath}
\usepackage{latexsym}
\usepackage[usenames]{color}

\numberwithin{equation}{section}

\newcommand{\be}{\begin{equation}}
\newcommand{\ee}{\end{equation}}
\newcommand{\beaa}{\begin{eqnarray*}}
\newcommand{\eeaa}{\end{eqnarray*}}
\newcommand{\bea}{\begin{eqnarray}}
\newcommand{\eea}{\end{eqnarray}}
\newcommand{\lbl}{\label}

\newcommand{\bei}{\begin{itemize}}
\newcommand{\eei}{\end{itemize}}

\newcommand{\bd}{\bold}

\newtheorem{theorem}{ \noindent T{\footnotesize HEOREM}}
\newtheorem{prop}{ \noindent P{\footnotesize ROPOSITION}}[section]
\newtheorem{lemma}{ \noindent L{\footnotesize EMMA}}[section]

\newcommand\td{\overset{d}{\to}}
\begin{document}

\title{Spectral Radii of Large Non-Hermitian Random Matrices}
\author{Tiefeng Jiang$^{1}$
and  Yongcheng Qi$^2$\\
University of Minnesota}

\date{}
\maketitle

\footnotetext[1]{School of Statistics, University of Minnesota, 224 Church
Street, S. E., MN55455, USA, jiang040@umn.edu. 
The research of Tiefeng Jiang was
supported in part by NSF Grant DMS-1209166 and DMS-1406279.}

\footnotetext[2]{Department of Mathematics and Statistics, University of Minnesota Duluth, MN 55812, USA, yqi@d.umn.edu.
The research of Yongcheng Qi was
supported in part by NSF Grant DMS-1005345.}

\begin{abstract}
\noindent By using the independence structure of points following a determinantal point process, we study the radii of the spherical ensemble, the truncation of the circular unitary ensemble and the product ensemble with parameter $n$ and $k$. The limiting distributions of the three radii are obtained. They are not the Tracy-Widom distribution. In particular, for the product ensemble, we show that the limiting distribution has a transition phenomenon: when $k/n\to 0$, $k/n\to \alpha\in (0,\infty)$ and $k/n\to \infty$, the liming distribution is the Gumbel distribution, a new distribution $\mu$ and the logarithmic normal distribution, respectively. The cumulative distribution function (cdf) of $\mu$ is the infinite product of some normal distribution functions. Another new distribution $\nu$ is also obtained for the spherical ensemble such that the cdf of $\nu$ is the infinite product of the cdfs of some Poisson-distributed random variables.  \\\\\\\\
\end{abstract}




\newpage

\section{Introduction}\lbl{Introduction}
\setcounter{equation}{0}
The largest eigenvalues of the three Hermitian matrices (Gaussian orthogonal ensemble, Gaussian unitary ensemble and Gaussian symplectic ensemble) are proved to converge to the Tracy-Widom laws by Tracy and Widom (1994, 1996).  Since then there have been very active research in this direction. For example,  Baik et al. (1999) establish a connection between the longest increasing subsequence problem and the Tracy-Widom law. The relationships among the largest eigenvalues, combinatorics, growth processes, random tilings and the determinantal point processes are found [see, e.g., Tracy-Widom (2002) and Johansson (2007) and the literature therein]. In the studies of the high-dimensional statistics,   Johnstone (2001, 2008) and Jiang (2009) prove that the largest eigenvalues of Wishart and Jacobi matrices converge to the Tracy-Widom law. Ram\'{\i}rez et al.  (2011) obtain the asymptotic distribution of the largest eigenvalues of beta-Hermite ensemble. Recently, a research interest is the universality of the largest eigenvalues of non-Gaussian matrices; see, for example, Tao and Vu (2011), Erd\H{o}s et al. (2012) and the references therein.

In this paper we will study the largest absolute values of the eigenvalues of some non-Hermitian matrices.  Initiated by Ginibre (1965) for the study of Gaussian random matrices (real, complex and symplectic), the interest has continued and theoretical results are found to have many applications in quantum chromodynamics, chaotic quantum systems and growth processes; see more descriptions from the paper by Akemann, Baik and Francesco  (2001). The applications also include dissipative quantum maps [Haake (2010)] and fractional quantum-Hall effect [Di Francesco et al.  (1994)]. We refer the readers to Khoruzhenko and  Sommers (2001) for more details.

For a matrix $\bd{M}$ with eigenvalues $z_1, \cdots, z_n$, the quantity $\max_{1\leq j \leq n}|z_j|$ is refereed to as the spectral radii of $\bd{M}$. In their pioneer work by  Rider (2003, 2004) and Rider and Sinclair (2014), the spectral radius of the real, complex and symplectic Ginibre ensembles are studied. For the complex Ginibre ensemble, it is shown that the spectral  radius converges to the Gumbel distribution. This phenomenon is very different from the Tracy-Widom distribution. The key observation is that the absolute values of the eigenvalues of the complex Ginibre ensemble are independent random variables with the Gamma distributions. The independence property is firstly observed by Kostlan (1992). Later it is found that the independence phenomenon is true not only  for the complex Ginibre ensemble, but also true for other complex-valued determinantal point processes; see, for example, Hough et al. (2009) for further details.

In this paper, we will study the largest radii of three rotation-invariant and non-Hermitian random matrices: the spherical ensemble $\bd{A}^{-1}\bd{B}$ where $\bd{A}$ and $\bd{B}$ are independent complex Ginibre ensembles, the truncation of circular unitary ensemble, and the product ensemble $\prod_{j=1}^k\bd{X}_j$ where $\bd{X}_1, \cdots,  \bd{X}_k$ are independent $n\times n$ complex Ginibre ensembles. The spectral  radii of the first one converges to a new distribution $\nu$, that of the second one converges to the Gumbel distribution, and that of the third one, depending on the ratio $\alpha:=\lim_{n\to\infty}k_n/n$, converges to the Gumbel distribution when $\alpha=0$, a new distribution $\mu$  when $\alpha\in (0, \infty)$ and the logarithmic normal distribution when $\alpha=\infty.$


Our analysis of the spectral radius is based on the following result. It  is a special case of Theorem 1.2 from  Chafa\"{\i} and P\'{e}ch\'{e} (2014) which is another version of  Theorem 4.7.1 from Hough et al. (2009).
\begin{lemma}\lbl{independence}(Independence of radius) Assume the density function of $(Z_1, \cdots, Z_n)\in \mathbb{C}^n$ is proportional to  $\prod_{1\leq j < k
\leq n}|z_j-z_k|^2\cdot \prod_{j=1}^n\varphi(|z_j|)$, where  $\varphi(x)\geq 0$ for all $x\geq 0.$ Let $Y_1, \cdots, Y_n$ be independent r.v.'s such that the density of $Y_j$ is proportional to  $y^{2j-1}\varphi(y)I(y\geq 0)$ for each $1\leq j \leq n.$ Then,  $g(|Z_1|, \cdots, |Z_n|)$ and $g(Y_1, \cdots, Y_n)$ have the same distribution for any symmetric function $g(y_1, \cdots, y_n)$.
\end{lemma}
Chafa\"{\i} and P\'{e}ch\'{e} (2014) also give two general results in their Theorems 1.3 and 1.4 to show the following: if the density function of the eigenvalues of a non-Hermitian random matrix is the same as that in Lemma \ref{independence}, under certain restrictions on $\varphi(x)$, the limiting distribution of the spectral radii is the Gumbel distribution. Their results do not apply to our three ensembles since our models do not meet their restrictions.

 Now we  present our results on the three ensembles in Subsections \ref{sub1}, \ref{sub2} and \ref{Pro_ensemble}, respectively. After this the strategy of the proofs and some comments are given.


\subsection{Spherical Ensemble}\lbl{sub1}
Let $\bd{A}$ and $\bd{B}$ be two $n\times n$ matrices and all of the $2n^2$ entries of the matrices are i.i.d. $\mathbb{C}N(0, 1)$-distributed random variables. Then $\bd{A}^{-1}\bd{B}$ is called a spherical ensemble [Hough et al. (2009)]. It has a connection to the matrix $F$ distribution in statistics literature; see, for instance, p. 331 from Eaton (2007).
Let $z_1, \cdots, z_n$ be the eigenvalues of $\bd{A}^{-1}\bd{B}$. Then their joint probability density function is given by
\bea\lbl{elastic}
C\cdot\prod_{j<k}|z_j-z_k|^2\cdot \prod_{k=1}^n\frac{1}{(1+|z_k|^2)^{n+1}}
\eea
where $C$ is a normalizing constant; see, for example,  Krishnapur (2009). The joint density of $z_1, \cdots, z_n$ of the real analogue of the spherical ensemble $\bd{A}^{-1}\bd{B}$, where $\bd{A}$ and $\bd{B}$ are i.i.d. real Ginibre ensembles, is given by Forrester and Nagao (2008) and Forrester and Mays (2011).

The empirical distribution of the eigenvalues has an asymptotic distribution $\mu$ with density $\frac{1}{\pi(1+|z|^2)^2}$ (Bordenave,  2011). When mapping the eigenvalues on the complex plane to the Riemann sphere  through the stereographic projection, the induced (pushforward) measure of $\mu$ is the uniform distribution on the sphere.  The spectral distribution of the singular values of $\bd{A}^{-1}\bd{B}$, which is the same as the eigenvalues of the $F$-matrix $(\bd{A}\bd{A}^*)^{-1}(\bd{B}\bd{B}^*)$, converges weakly to a non-random distribution; see, for instance, Wachter (1980) and  Bai et al. (1987).

In this paper, we say
\beaa
X_n\ \ \mbox{converges weakly to the cdf } F(x) \mbox{ or a random variable } X
\eeaa
if the probability distribution of $X_n$ converges weakly to that generated by the cumulative distribution function (cdf) $F(x)$ or $X.$  Now we study the spectral radius.
\begin{theorem}\lbl{Theorem_Spherical} Let $z_1, \cdots, z_n$ have the density as in (\ref{elastic}).
Define $H_k(x)=e^{-x}\sum^{k-1}_{j=0}\frac{x^j}{j!}$ for $k\ge 1$.  Then $\frac{1}{\sqrt{n}}\max_{1\leq j \leq n}|z_j|$ converges weakly to  probability distribution  function $H(x)=\prod^\infty_{k=1}H_k(x^{-2})$ for $x>0$ and $H(x)=0$ for $x \leq 0.$
\end{theorem}

Observe that $H_k(x)$ is the cdf $P(\mbox{Poi}(x) \leq k-1)$ for each $k\geq 1$, where $\mbox{Poi}(x)$ is a  Poisson random variable with parameter $x>0.$ So $H(x)$ is the product of those cdfs evaluated at $x^{-2}$.

Johnstone (2008) proves that, under a trivial transformation, the largest singular value of the $F$-matrix $(\bd{A}\bd{A}^*)^{-1}(\bd{B}\bd{B}^*)$ asymptotically follows a Tracy-Widom distribution. Here, the spectral radius converges weakly to the new distribution $H(x)$.  Now we examine the tail probability of the distribution function $H(x)$ as in Theorem \ref{Theorem_Spherical}. In fact, we have
\bea\lbl{smartyy}
1-H(x)\sim \frac{1}{x^{2}}
\eea
as $x\to +\infty$. So $H(x)$ is heavy-tailed. This property will be verified in Section \ref{Flyingcolor}.

\subsection{Truncation of Circular Unitary Ensemble}\lbl{sub2}

Now we consider the truncation of the circular unitary ensemble. Let $\bd{U}$ be an $n\times n$ Haar-invariant unitary matrix [see, e.g.,  Diaconis and Evans (2001) and Jiang (2009, 2010)]. For $n>p\geq 1$, write
\begin{eqnarray*}
\bd{U}=\begin{pmatrix}
\bd{A}\ \ \bd{C}^*\\
\bd{B}\ \ \bd{D}
\end{pmatrix}
\end{eqnarray*}
where $\bd{A}$, as a truncation of $\bd{U}$, is a $p\times p$ submatrix. Let $z_1, \cdots, z_p$ be the eigenvalues of $\bd{A}$. It is known from Zyczkowski and Sommers (2000) that their density function is
\bea\lbl{Good_forever}
C\cdot\prod_{1\leq j<k\leq p}|z_j - z_k|^2\prod_{j=1}^{p}(1-|z_j|^2)^{n-p-1}
\eea
where $C$ is a normalizing constant. Assuming $c=\lim\frac{p}{n},$ \.{Z}yczkowski and Sommers (2000) show that the empirical distribution of $z_i$'s converges to the distribution with density proportional to $\frac{1}{(1-|z|^2)^2}$ for $|z|\leq c$ if $c\in (0,1).$ Dong et al. (2012) prove that the empirical distribution goes to the circular law and the arc law as $c=0$ and $c=1$, respectively.

Collins (2005) proves that $\bd{A}^*\bd{A}$ forms a Jacobi ensemble. Johansson (2000) and Jiang (2009) show that a transform of the largest eigenvalue of $\bd{A}^*\bd{A}$ converges weakly to the Tracy-Widom distribution. For the spectral radius $\max_{1\leq j \leq p}|z_j|$ of $\bd{A}$ itself, we obtain the following result.


\begin{theorem}\label{Theorem_hyperbolic}
Assume that $z_1, \cdots, z_p$ have density as in (\ref{Good_forever}) and there exist constants $h_1, h_2\in (0,1)$ such that $h_1< \frac{p}{n}< h_2$ for all $n\geq 2.$ Then $(\max_{1\leq j \leq p}|z_j|-A_n)/B_n$ converges weakly to the cdf $\Lambda(x)=\exp(-e^{-x})$, $x\in \mathbb{R}$,
where $A_n=c_n+\frac{1}{2}(1-c_n^2)^{1/2}(n-1)^{-1/2}a_n$,  $B_n=\frac{1}{2}(1-c_n^2)^{1/2}(n-1)^{-1/2}b_n$,
\[
c_n=\Big(\frac{p-1}{n-1}\Big)^{1/2},~~ b_n=b\Big(\frac{nc_n^2}{1-c_n^2}\Big), ~~a_n=a\Big(\frac{nc_n^2}{1-c_n^2}\Big)
\]
with
\begin{equation}\label{abx}
a(y)=(\log y)^{1/2}-(\log y)^{-1/2}\log(\sqrt{2\pi}\log )\ \mbox{ and } \ b(y)=(\log y)^{-1/2}
\end{equation}
for $y>3$.
\end{theorem}
Trivially, in the above theorem, $\{A_n;n \geq 3\}$ is bounded and $B_n$ has the scale of $(n\log n)^{-1/2}$.

\subsection{Product Ensemble}\lbl{Pro_ensemble}

Given integer $k\geq 1$. Assume $\bd{X}_1, \cdots, \bd{X}_k$ are i.i.d. $n\times n$ random matrices and the $n^2$ entries of $\bd{X}_1$ are i.i.d. with distribution $\mathbb{C}N(0, 1).$ Let $z_1, \cdots, z_n$ be the eigenvalues of the product $\prod_{j=1}^k\bd{X}_j$. It is known that their joint density function is
\bea\label{up_stairs}
C\prod_{1\leq j <l \leq n}|z_j - z_l|^2\prod_{j=1}^nw_{k}(|z_j|)
\eea
where $C$ is a normalizing constant and $w_k(z)$  is
given by the Meijer G-function
with
\bea\lbl{smell_good}
w_k(z)=\pi^{k-1}G_{n, 0}^{0,  n}\Big(\begin{matrix}- \\ {\overrightarrow{\bd{0}}} \end{matrix}\Big||z|^2\Big).
\eea
This formula seems not easy to understand at the first sight. However, the function admits an easily recursive formula $w_1(z)=\exp(-|z|^2)$ and
\begin{equation}\label{wk}
w_{k}(z)=2\pi\int^\infty_0w_{k-1}\Big(\frac{z}{r}\Big)\exp(-r^2)\frac{dr}{r}
\end{equation}
for all integer $k\geq 2$; see, for example, Akemann and Burda (2012).

A paper in Arxiv by G\"{o}tze and Tikhomirov says  that the empirical distribution of $z_j/n^{k/2},\, 1\leq j \leq n$, in the sense of mean value, converges to a distribution with density $\frac{1}{k\pi}|z|^{\frac{2}{k}-2}$ for $|z|\leq 1.$ Later, Bordenave (2011), O'Rourke and Soshnikov (2011) and O'Rourke et al. (2014) further generalize this result to the almost sure convergence. The Gaussian case was first considered by Burda et al. (2010) and  Burda (2013) through investigating the limit of the kernel of a determinantal point process. Second, the empirical distribution of the singular values of $\bd{X}_1\bd{X}_2/n$ converges weakly to a non-random distribution, see, for instance, Theorem 2.10 from Bai (1999).

Now we consider the largest radius and the result is given below.  We allow $k$ changes with $n$ in this paper. First, we need some notation. Let $\Phi$ denote the cumulative distribution function of $N(0, 1)$. For $\alpha\in (0, \infty)$,
define
\[
\Phi_\alpha(x)=\prod^\infty_{j=0}\Phi\Big(x+j\alpha^{1/2}\Big),
\]
and $\Phi_\infty(x)=\Phi(x)$. The digamma function $\psi$ is defined by
\bea\lbl{maomi}
\psi(z)=\frac{d}{dz}\log\Gamma(z)=\frac{\Gamma'(z)}{\Gamma(z)}
\eea
where $\Gamma(z)$ is the Gamma function.

\begin{theorem}\lbl{Theorem_product} Let $k=k_n$ be a sequence of positive integers. The following holds.\\
\noindent{(a).}  If $\lim_{n\to\infty}k_n/n=0$, particularly for $k_n\equiv k$, then $\alpha_n\big(n^{-k_n/2}\max_{1\le j\le n}|z_j|-1\big)-\beta_n$ converges weakly to the cdf  $\exp(-e^{-x})$,
where
\beaa
\alpha_n= \Big(\frac{n}{k_n}\log \frac{n}{k_n}\Big)^{1/2}\ \ \ \mbox{and}\ \ \ \beta_n=\log \frac{n}{k_n}-\log \log \frac{n}{k_n}-\frac{1}{2}\log (2\pi).
\eeaa

\noindent{(b).} If $\lim_{n\to\infty}k_n/n=\alpha\in (0,\infty)$, then
\beaa
\frac{\max_{1\le j\le n}|z_j|}{n^{k_n/2}}\ \ \mbox{converges weakly to the cdf}\ \  \Phi_{\alpha}\Big(\frac{1}{2}\alpha^{1/2}+2\alpha^{-1/2}\log x\Big),~~~~x>0.
\eeaa

\noindent{(c).} If $\lim_{n\to\infty}k_n/n=\infty$, then
\beaa
\frac{\max_{1\le j\le n}\log |z_j|-k_n\psi(n)/2}{\sqrt{k_n/n}/2}\ \ \mbox{converges weakly to}\  N(0, 1).
\eeaa
\end{theorem}

Taking $k=1$ in (a) of Theorem \ref{Theorem_product}, the corresponding limiting result is obtained by Rider (2003). Here we not only get the result for finite $k$, but for all possible range of $k_n$, which leads to the three transition zones: $k_n/n\to \alpha$ with $\alpha=0$, $\alpha\in (0, \infty)$ and $\alpha=\infty.$

As mentioned below Lemma \ref{independence}, Theorems 1.3 and 1.4 from Chafa\"{\i} and P\'{e}ch\'{e} (2014) conclude that the limiting distributions are always the Gumbel. The two theorems do not imply any of our results. Although the  limiting distributions in Theorem \ref{Theorem_hyperbolic} and  case (a) of Theorem \ref{Theorem_product} are the Gumbel, since the density functions in (\ref{Good_forever}) and (\ref{up_stairs}) have two parameters $k$ and $n$ with $k$ depending on $n$,  their assumptions are not satisfied.

Second, let $\eta_{\max}(\bd{X}_i)$ be the largest singular value of $\bd{X}_i$ for each $i$. It is proved that $\eta_{\max}(\bd{X}_1)/\sqrt{n}$ $\to 2$ in probability (see, e.g., Bai, 1999). Let $\lambda_1(\bd{X}_i), \cdots, \lambda_n(\bd{X}_i)$ be the eigenvalues of $\bd{X}_i$ for each $i$. Obviously, $\max_{1\leq j \leq n}|\lambda_j(\bd{X}_i)| \leq \eta_{\max}(\bd{X}_i)$ for each $i$. From Rider (2003) or (a) of Theorem \ref{Theorem_product}, we know $\max_{1\leq j \leq n}|\lambda_j(\bd{X}_i)|/\sqrt{n} \to 1$ in probability for each $i$. It is interesting to see the first limit is $2$ and the second is $1$. Now, the assertion (b) of Theorem \ref{Theorem_product} says that, though $\max_{1\leq j \leq n}|\lambda_j(\bd{X}_i)|/\sqrt{n} \to 1$ for each $i$, the radius of $\prod_{i=1}^k(\bd{X}_i/\sqrt{n}\,)$ goes to a distribution with support $[0, \infty).$

Finally, let us look at the tail behavior of the distribution  in (b) of Theorem \ref{Theorem_product}. In fact we have
\bea\lbl{wisconsin}
1-\Phi_{\alpha}\Big(\frac{1}{2}\alpha^{1/2}+2\alpha^{-1/2}\log y\Big)&\sim& \frac{C}{x\log x}e^{-2(\log x)^2/\alpha}
\eea
as $x\to + \infty$, where $C=\frac{\sqrt{\alpha} e^{-\alpha/8}}{2\sqrt{2\pi}}$. It is different from that of $e^{N(0,1)}$, the standard logarithmic normal distribution: $P(e^{N(0,1)}\geq x)\sim \frac{1}{\sqrt{2\pi}\,\log x}e^{-(\log x)^2/2}$ as $x\to +\infty.$ This will be verified in Section \ref{Flyingcolor}.\\

\noindent\textbf{Strategy of the proofs}. By using Lemma \ref{independence}, the absolute values of eigenvalues $|z_i|$'s are ``independent". So we are dealing with the maxima of independent random variables with different distributions. The first step is to identify the distribution of each random variable. For example, for the product ensemble in Section \ref{Pro_ensemble}, $|z_j|$ has the same distribution as the product of some i.i.d. random variables with Gamma distributions (Lemma \ref{Catroll}). Then we analyze  the  tail probabilities of the product of random variables carefully through moderate deviations (Proposition \ref{prop5}). This step costs the major effort.  \\

\noindent\textbf{Comments}:

1. It is noteworthy to mention that, though the main idea is analyzing the maxima of independent random variables, the proofs are not trivial. In the classical study of the maxima of i.i.d. random variables, the limiting distributions are only of three types:  Fr\'{e}chet distribution, Gumbel distribution and Weibull distribution; see, for example, Resnick (2007). However, the limiting distributions appeared in Theorems \ref{Theorem_Spherical} and (b) of Theorem \ref{Theorem_product} are new.

2. The eigenvalues of the three random matrices investigated in this paper are rotation-invariant. This special property gives us the advantage of independence by Lemma \ref{independence}. When the eigenvalues are not of the invariant property, it seems there have no good understanding on  the largest radii. For example, if $z_1, \cdots, z_n$ have joint density
\beaa
f(z_1, \cdots, z_n) =C\cdot \prod_{1\leq j < k \leq n}|z_j-z_k|^2\cdot \exp\Big\{-n\sum_{j=1}^n\Big(\frac{(\mbox{Re} z_j)^2}{1+\tau} - \frac{(\mbox{Im} z_j)^2}{1-\tau}\Big) \Big\}
\eeaa
where $\tau \in (-1, 1)$ is a parameter and $C$ is a normalizing constant [Lemma 4 from Petz and Hiai (1998)]. See also a similar example on page 3403 from Rider (2003) and  (1.1) from the Arxiv paper by Kuijlaars and L\'{o}pez.

3. In this paper, we work on matrices with complex Gaussian entries. A similar study may be done for matrices with real and symplectic Gaussian random variables. For example, we know from Ginibre (1965), Lehmann and Sommers (1991) and  Edelman (1997) that the densities of the eigenvalues of real and symplectic Ginibre ensembles are also explicit. Rider (2003) and Rider and Sinclair (2014) obtain the limiting distributions of the largest radii for the real and symplectic cases. It is possible that our current work can be carried out to the three real and symplectic analogues: the spherical ensemble $\bd{A}^{-1}\bd{B}$ where $\bd{A}$ and $\bd{B}$ are real or symplectic Ginibre ensembles [further information can be seen from Forrester and Nagao (2008) and Forrester and Mays (2011)], the truncation of Haar-invariant orthogonal or symplectic matrices [see, e.g.,  Jiang (2010)] and $\prod_{j=1}^k\bd{X}_j$ where $\bd{X}_1, \cdots, \bd{X}_k$ are i.i.d. real or symplectic Ginibre ensembles.

4. Tracy and Widom (1994, 1996) prove that the largest eigenvalues of the Gaussian orthogonal, unitary and symplectic ensembles converge to the Tracy-Widom laws. Recently there have been an active research on the universality of the eigenvalues of non-Gaussian matrices; see, for example, Tao and Vu (2011), Erd\H{o}s et al. (2012) and the references therein.  In particular, Erd\H{o}s et al. generalize the results by Tracy-Widom to the matrices with  non-Gaussian entries. Our Theorems \ref{Theorem_Spherical}, \ref{Theorem_hyperbolic} and \ref{Theorem_product} consider the eigenvalues of matrices with Gaussian entries. It will be interesting to study the universality of the three results for the matrices with non-Gaussian entries.

Finally, the organization of the rest of paper is as follows. We will prove  Theorems \ref{Theorem_Spherical}, \ref{Theorem_hyperbolic} and \ref{Theorem_product} in Sections  \ref{wo}, \ref{nimen} and \ref{jinren}, respectively. The verifications of (\ref{smartyy}) and (\ref{wisconsin}) are given in Section \ref{Flyingcolor}.





\section{Proofs}
In this section, we will prove Theorems \ref{Theorem_Spherical}, \ref{Theorem_hyperbolic} and \ref{Theorem_product} in each subsection.

\subsection{The Proof of Theorem \ref{Theorem_Spherical}}\lbl{wo}

\noindent We start with a lemma.
\begin{lemma}\label{lem1} Let $a_{ni}\in [0,1)$ be constants for $i\ge 1, n\ge 1$ and $\sup_{n\ge 1, i\ge 1}a_{ni}<1$.
For each $i\ge 1$,  $a_i:=\lim_{n\to\infty}a_{ni}$.  Assume $c_n:=\sum_{i=1}^{\infty}a_{ni}<\infty$ for each $n\ge 1$ and
$c:=\sum_{i=1}^{\infty}a_i<\infty$, and $\lim_{n\to\infty}c_n=c$.  Then
\begin{equation}\label{prod}
\lim_{n\to\infty}\prod^\infty_{i=1}(1-a_{ni})=\prod^\infty_{i=1}(1-a_i).
\end{equation}
\end{lemma}

\noindent\textbf{Proof.}  Note that $\prod^\infty_{i=1}(1-a_{ni})$ and $\prod^\infty_{i=1}(1-a_i)$ are well defined, and
$\prod^\infty_{i=1}(1-a_{ni})>0$ for each $n\ge 1$ and $\prod^\infty_{i=1}(1-a_i)>0$.  It suffices to show that
\begin{equation}\label{2sums}
\lim_{n\to\infty}\sum^\infty_{i=1}\log(1-a_{ni})=\sum^\infty_{i=1}\log(1-a_i).
\end{equation}

Note that for each fixed $k>1$
\begin{eqnarray*}
\limsup_{n\to\infty}\sum^\infty_{i=1}|a_{ni}-a_i|&\le &\limsup_{n\to\infty}\sum^k_{i=1}|a_{ni}-a_i|
       +\limsup_{n\to\infty}\sum^\infty_{i=k+1}|a_{ni}-a_i|\\
       &=&\limsup_{n\to\infty}\sum^\infty_{i=k+1}|a_{ni}-a_i|\\
       &\le &\limsup_{n\to\infty}\sum^\infty_{i=k+1}a_{ni}
       +\sum^\infty_{i=k+1}a_i\\
        &=&\limsup_{n\to\infty}(c_n-\sum^k_{i=1}a_{ni})
       +c-\sum^k_{i=1}a_i\\
              &=&2(c-\sum^k_{i=1}a_i),
\end{eqnarray*}
which goes to zero as $k\to\infty$. Therefore,
we have
\[
\lim_{n\to\infty}\sum^\infty_{i=1}|a_{ni}-a_i|=0.
\]
Set $a=\sup_{n\ge 1, i\ge 1}a_{ni}$. Then $0\leq a<1.$ It follows that
\begin{eqnarray*}
|\sum^\infty_{i=1}\log(1-a_{ni})-\sum^\infty_{i=1}\log(1-a_i)|
&\le &\sum^\infty_{i=1}|\log(1-a_{ni})-\log(1-a_i)|\\
&= &\sum^\infty_{i=1}|\int^{a_{ni}}_{a_i}\frac1{1-t}dt|\\
&\le& \frac{1}{1-a}\sum^\infty_{i=1}|a_{ni}-a_i| \to 0
\end{eqnarray*}
as $n\to\infty$, proving \eqref{2sums}.  \hfill$\Box$\\

\noindent\textbf{Proof of Theorem \ref{Theorem_Spherical}}. By Lemma \ref{independence} and (\ref{elastic}), $\max_{1\leq j \leq n}|z_j|$ and $\max_{1\leq j \leq n}Y_{nj}$ have the same distribution, where $Y_{n1}, \cdots, Y_{nn}$ are independent such that $Y_{nj}$ has the probability density function (pdf) proportional to $y^{2j-1}(1+y^2)^{-(n+1)}I(y\geq 0)$ for $1\leq j \leq n$. Thus, to prove the theorem, it suffices to show
\bea\lbl{dong_beauty}
\lim_{n\to\infty}P\Big(\frac{1}{\sqrt{n}}\max_{1\leq j \leq n}Y_{nj} \leq x\Big)=H(x)
\eea
for each $x>0.$

Let $X_i$,  $i\ge 1$ be a sequence of i.i.d. random variables with cumulative distribution function (cdf) $F$.
Let $X_{1:n}\le X_{2:n}\le\cdots\le X_{n:n}$ be the order statistics of $X_1, X_2, \cdots, X_n$ for each $n\ge 1$.  Then
from page 14 on the book by Balakrishnan and Cohen (1991), we know that the cdf of $X_{i:n}$ is given by
\begin{equation}\label{cdf}
F_{i:n}(x)=\sum^n_{r=i}{{n}\choose{r}}F(x)^r(1-F(x))^{n-r}=\frac{n!}{(i-1)!(n-i)!}\int^{F(x)}_0
t^{i-1}(1-t)^{n-i}dt
\end{equation}
 for each $1\le i\le n$.
 If $F$ has a probability density function  $f$, then the pdf of $X_{i:n}$ is given by
  \begin{equation}\label{pdf}
f_{i:n}(x)=\frac{n!}{(i-1)!(n-i)!}
F(x)^{i-1}(1-F(x))^{n-i}f(x).
\end{equation}

The monotonicity of the order statistics implies that $F_{i:n}(x)$ is non-increasing in $i$ for each $x$, that is,
\begin{equation}\label{monotone}
F_{1:n}(x)\ge F_{2:n}(x)\ge\cdots\ge F_{n:n}(x).
\end{equation}

Let $\{u_n, ~n\ge 1\}$ be a sequence of constants such that $\lim_{n\to\infty}n(1-F(u_n))=:\tau\in (0,\infty)$. Write $\tau_n=n(1-F(u_n)).$ Then it follows
from the first equality in equation \eqref{cdf} that
\bea\label{fun}
F_{n-i+1: n}(u_n) &=& \sum_{r=n-i+1}^n\binom{n}{r}F(u_n)^r(1-F(u_n))^{n-r} \nonumber\\
& = & \sum_{j=0}^{i-1}\binom{n}{j}\Big(1-\frac{\tau_n}{n}\Big)^{n-j}\Big(\frac{\tau_n}{n}\Big)^{j} \nonumber\\
& = & \sum_{j=0}^{i-1}\frac{1}{j!}\cdot\prod_{l=1}^{j-1}\Big(1-\frac{l}{n}\Big)\cdot \Big(1-\frac{\tau_n}{n}\Big)^{n-j}(\tau_n)^{j} \nonumber\\
&\to& e^{-\tau}\sum_{j=0}^{i-1}\frac{\tau^j}{j!}= H_i(\tau)
\eea
as $n\to\infty$ for each fixed integer $i\ge 1$.

Now, we take $F(y)=\frac{y^2}{1+y^2}$ for $y>0$.  Fix $x>0$,
set $u_n=u_n(x)=\sqrt{n}x$. Then $\lim_{n\to\infty}n(1-F(u_n))=x^{-2}$. Then from \eqref{fun}
\begin{equation}\label{morefun}
\lim_{n\to\infty}F_{n-i+1: n}(u_n)=H_i(x^{-2})
\end{equation}
for each fixed integer $i\ge 1$. For each $n\ge 1$, define
\[
a_{ni}=\left\{
         \begin{array}{ll}
           1-F_{n-i+1:n}(u_n), & \hbox{if $1\le i\le n$;} \\
           0, & \hbox{if $i>n$.}
         \end{array}
       \right.
\]
Then it follows from \eqref{monotone} that $\sup_{n\ge 1, i\ge 1}a_{ni}=\sup_{n\ge 1}a_{n1}$. By the first identity in (\ref{cdf}),
\begin{eqnarray*}
a_{n1}&=&1-\Big(\frac{nx^2}{1+nx^2}\Big)^n=1-\Big(1+\frac{x^{-2}}{n}\Big)^{-n}
\end{eqnarray*}
is increasing in $n$.
Hence,
\[
\sup_{n\ge 1, i\ge 1}a_{ni}=\lim_{n\to \infty}a_{n1}=1-\exp(-x^{-2})\in (0,1).
\]
From \eqref{morefun} we have $\lim_{n\to\infty}a_{ni}=1-H_i(x^{-2})=:a_i$ for each $i\ge 1$.
Moreover, we have
\beaa
\sum^{\infty}_{i=1}a_{ni}&=& E\Big[\sum^n_{i=1}I(X_{n-i+1:n}>u_n)\Big]\\
& = & E\Big[\sum^n_{i=1}I(X_i>u_n)\Big]\\
&=& n(1-F(u_n))\to x^{-2}
\eeaa
and
\[
\sum^{\infty}_{i=1}a_{i}=\sum^{\infty}_{i=1}(1-H_i(x^{-2}))=\exp(-x^{-2})\sum^{\infty}_{i=1}\sum^\infty_{k= i}\frac{(x^{-2})^k}{k!}.
\]
By exchanging the ordering of the sums, we know $\sum^{\infty}_{i=1}\sum^\infty_{k= i}\frac{(x^{-2})^k}{k!}=\sum^{\infty}_{k=1}\sum^k_{i=1}\frac{(x^{-2})^k}{k!}=\sum^{\infty}_{k=1}\frac{(x^{-2})^k}{(k-1)!}=x^{-2}\exp\{x^{-2}\}.$
It follows that $\sum^{\infty}_{i=1}a_{i}=x^{-2}.$
By Lemma~\ref{lem1}, we have
\[
\lim_{n\to\infty}\prod^n_{i=1}F_{i:n}(\sqrt{n}x)=\lim_{n\to\infty}\prod^\infty_{i=1}(1-a_{ni})
=\prod^\infty_{i=1}H_i(x^{-2})=H(x)~~~\mbox{ for }x>0.
\]
From \eqref{pdf} we obtain the pdf of $X_{j:n}$ given by
\[
f_{j:n}(y)=\frac{2n!}{(j-1)!(n-j)!}y^{2j-1}(1+y^2)^{-(n+1)},  ~~~y>0,
\]
which is also the pdf of $Y_{nj}$. Therefore, we have

\begin{eqnarray*}
 \lim_{n\to\infty}P\Big(\frac{1}{\sqrt{n}}\max_{1\leq j \leq n}Y_{nj} \leq x\Big)
&=& \lim_{n\to\infty}\prod^n_{j=1}P(Y_{nj}\le \sqrt{n}x)\\
&=&\lim_{n\to\infty}\prod^n_{j=1}F_{j:n}(\sqrt{n}x)=H(x)
\end{eqnarray*}
for $x>0.$ This completes the proof of Theorem \ref{Theorem_Spherical}. \hfill $\blacksquare$

\subsection{The Proof of Theorem \ref{Theorem_hyperbolic}}\lbl{nimen}

Notation:  $C_n\sim D_n$ as $n\to\infty$ implies $\lim_{n\to\infty}\frac{C_n}{D_n}=1$;  $C_n(t)\sim D_n(t)$ uniformly over $t\in
T_n$ implies $\lim_{n\to\infty}\sup_{t\in T_n }|\frac{C_n(t)}{D_n(t)}-1|=0$;  $C_n(t)=O(D_n(t))$ uniformly over $t\in T_n$
implies $\sup_{t\in T_n}|\frac{C_n(t)}{D_n(t)}|$ is bounded;  $C_n(t)=o(D_n(t))$ uniformly over $t\in T_n$
implies $\sup_{t\in T_n}|\frac{C_n(t)}{D_n(t)}|$ converges to zero as $n\to\infty$.

For random variables $\{X_n;\, n\geq 1\}$ and constants $\{a_n;\, n\geq 1\}$, we write $X_n=O_P(a_n)$ if $\lim_{x\to +\infty}\lim_{n\to\infty}P(|\frac{X_n}{a_n}|\geq x)=0$. In particular, if $X_n=O_P(a_n)$ and $\{b_n;\, n\geq 1\}$ is a sequence of constants with $\lim_{n\to\infty}b_n=\infty$, then $\frac{X_n}{a_nb_n}\to 0$ in probability as $n\to\infty$.

Recall $\phi(x)=\frac{1}{\sqrt{2\pi}}e^{-x^2/2}$  and $\Phi(x)=\frac{1}{\sqrt{2\pi}}\int_{-\infty}^xe^{-t^2/2}\,dt$ for $x\in \mathbb{R}.$ Let also $a(x)$ and  $b(x)$ be as in Theorem \ref{Theorem_hyperbolic}.

\begin{lemma}\label{seed}
Let $\{j_n,  n\ge 1\}$ and $\{x_n, n\geq 1\}$ be positive numbers with
$\lim_{n\to\infty}x_n=\infty$
and  $\lim_{n\to\infty}j_n x_n^{-1/2}(\log x_n)^{1/2}=\infty$. Let
$\{c_{n,j}, 1\le j\le j_n, n\ge 1\}$ be real numbers
with $\lim_{n\to\infty}\max_{1\le j\le j_n}|c_{n,j}x_n^{1/2}-1|=0$. Then, for all $y\in \mathbb{R}$,
\bea
& & \lim_{n\to\infty}\sum_{j=1}^{j_n}\big(1-\Phi((j-1)c_{n,j}+a(x_n)+b(x_n)y\big)=e^{-y};\label{normal-app1}\\
& & \lim_{n\to\infty}\sum_{j=1}^{j_n}\frac{1}{(j-1)c_{n,j}+a(x_n)+b(x_n)y}\phi((j-1)c_{n,j}+a(x_n)+b(x_n)y\big)=e^{-y}.
\label{normal-app2}  \ \ \ \ \ \ \ \ \ \
\eea
\end{lemma}

\noindent\textbf{Proof.} From definition, it is easy to see that $\lim_{n\to\infty}a(x_n)= +\infty$ and $\lim_{n\to\infty}b(x_n)= 0$ and $\min_{1\leq j \leq j_n}c_{n,j}>0$ as $n$ is large enough. Thus, $\min_{1\leq j \leq j_n}[(j-1)c_{n,j}+a(x_n)+b(x_n)y]\to +\infty$ as $n\to\infty$.  It is well known that  $1-\Phi(x)\sim \frac{\phi(x)}{x}$  as  $x\to\infty.$ Therefore, \eqref{normal-app1} follows from  \eqref{normal-app2}. Now let us prove \eqref{normal-app2}.

It is easy to verify that for $y\in \mathbb R$
\begin{equation}\label{aby}
\exp\big(-\frac{1}{2}(a(x_n)+b(x_n)y)^2\big)\sim \frac{\sqrt{2\pi}\log x_n}{\sqrt{x_n}}e^{-y}
\end{equation}
as $n\to\infty$. For large $n$, define
\[
l_n=\mbox{ the integer part of }\min\Big\{\frac{j_n^{1/2}x_n^{1/4}}{(\log x_n)^{1/4}}, \frac{x_n^{1/2}}{(\log x_n)^{1/4}}\Big\}.
\]
Then, as $n\to\infty$,
\begin{equation}\label{onln}
j_n>l_n\to \infty, ~~~ l_n x_n^{-1/2}(\log x_n)^{1/2}\to \infty~ \mbox{ and } ~ l_nx_n^{-1/2}\to 0.
\end{equation}
Fix $y\in \mathbb{R}$ and set
\[
u_{nj}=(j-1)c_{n,j}+a(x_n)+b(x_n)y, ~~~~1\le j\le j_n.
\]
Then we conclude the following facts:

\noindent{\bf Fact 1:} Uniformly over $1\le j\le l_n$,
\begin{equation}\label{1toln-u}
u_{nj}=(\log x_n)^{1/2}(1+o(1))
\end{equation}
by using the third assertion in (\ref{onln}) and
\begin{equation}\label{1toln-u2}
u_{nj}^2=(a(x_n)+b(x_n)y)^2+2(j-1)\frac{(\log x_n)^{1/2}}{{x_n}^{1/2}}(1+o(1))+o(1);
\end{equation}

\noindent{\bf Fact 2:} Uniformly over $l_n < j\le j_n$, which is different from the assumption on (\ref{1toln-u}) and (\ref{1toln-u2}),
\begin{equation}\label{lntojn-u}
u_{nj}\ge (\log x_n)^{1/2}(1+o(1))
\end{equation}
and
\begin{equation}\label{lntojn-u2}
u_{nj}^2\ge (a(x_n)+b(x_n)y)^2+2(j-1)\frac{(\log x_n)^{1/2}}{{x_n}^{1/2}}(1+o(1))+o(1).
\end{equation}
It then follows from \eqref{1toln-u}, \eqref{1toln-u2} and  \eqref{aby} that
\begin{eqnarray}
\sum^{l_n}_{j=1}\frac{1}{u_{nj}}\phi(u_{nj})
&\sim&\frac{\exp\big(-\frac{1}{2}(a(x_n)+b(x_n)y)^2\big)}{(2\pi\log x_n)^{1/2}}\sum^{l_n}_{j=1}
         \exp\big(-(j-1)\frac{(\log x_n)^{1/2}}{{x_n}^{1/2}}(1+o(1))\big) \nonumber\\
&\sim&\frac{\exp\big(-\frac{1}{2}(a(x_n)+b(x_n)y)^2\big)}{(2\pi\log x_n)^{1/2}}
\frac{1}{1-\exp\big(-\frac{(\log x_n)^{1/2}}{{x_n}^{1/2}}(1+o(1))\big)} \lbl{silver_bar}\\
&\sim&\frac{\exp\big(-\frac{1}{2}(a(x_n)+b(x_n)y)^2\big)}{(2\pi\log x_n)^{1/2}}
\frac{{x_n}^{1/2}}{(\log x_n)^{1/2}}
\sim e^{-y},\lbl{Gold_bar}
\end{eqnarray}
where the middle limit in (\ref{onln}) is used in the second step.
Similarly, it follows from \eqref{lntojn-u}, \eqref{lntojn-u2} and \eqref{aby} that
\begin{eqnarray*}
& & \sum^{j_n}_{j=l_n+1}\frac{1}{u_{nj}}\phi(u_{nj})\\
&\le &\frac{\exp\big(-\frac{1}{2}(a(x_n)+b(x_n)y)^2\big)}{(2\pi\log x_n)^{1/2}}\sum^{\infty}_{j=l_n+1}
         \exp\big(-(j-1)\frac{(\log x_n)^{1/2}}{{x_n}^{1/2}}(1+o(1))\big)\\
&\le &
\frac{\exp\big(-\frac{1}{2}(a(x_n)+b(x_n)y)^2\big)}{(2\pi\log x_n)^{1/2}}\frac{1}{1-\exp\big(-\frac{(\log x_n)^{1/2}}{{x_n}^{1/2}}(1+o(1))\big)}\\
&&~~~~~~~~~~~~~~~~~~~~~~~~~~~~~~~~~~~~~~~~~~~~~\times\exp\Big(-l_n\frac{(\log x_n)^{1/2}}{{x_n}^{1/2}}(1+o(1))\Big)\\
&=&O\Big(\exp\Big(-l_n\frac{(\log x_n)^{1/2}}{{x_n}^{1/2}}(1+o(1))\Big)\Big)\to 0
\end{eqnarray*}
by using \eqref{silver_bar} and \eqref{Gold_bar} in the equality and  the middle assertion in (\ref{onln}) in the last step. By adding up the above eqaution and (\ref{Gold_bar}), we obtain \eqref{normal-app2}.   \hfill $\blacksquare$\\

Let $\{U_i;\,  i\ge 1\}$ be a sequence of i.i.d. random variables uniformly distributed over $(0,1)$, and let $U_{1:n}\le U_{2:n}\le \cdots\le U_{n:n}$ be the order statistics of $U_1, U_2, \cdots, U_n$ for each $n\ge 1$. Recall $\cal{B}$ is the collection of all Borel sets on $\mathbb{R}$. The following lemma is a special case of Proposition 2.10 from  Reiss (1981).
\begin{lemma}\label{reiss}  There exists a constant $C>0$ such that for all $r>k\geq 1$,
\beaa
& & \sup_{B\in \cal{B}}\Big|P\Big(\frac{r^{3/2}}{\sqrt{(r-k)k}}\big(U_{r-k+1:r}-\frac{r-k}{r}\big)\in B\Big)-\int_B(1+l_1(t)+l_2(t))\phi(t)dt \Big|\\
&\le &
C\cdot\Big(\frac{r}{(r-k)k}\Big)^{3/2}
\eeaa
where
for $i=1,2$, $l_i(t)$ is a polynomial in $t$ of degree $\le 3i$, depending on $r$ and $k$,  and all of its coefficients  are of order $O(\big(\frac{r}{(r-k)k}\big)^{i/2})$.
\end{lemma}

\noindent\textbf{Proof of Theorem \ref{Theorem_hyperbolic}}. Review the density formula in (\ref{Good_forever}). Set $m_n=n-p.$ For ease of notation, we sometimes write $m$ for $m_n.$ By assumption, $h_2'< \frac{m_n}{n}< h_1'$ for all $n\geq 2$ where $h_i'=1-h_i\in (0,1)$ for $i=1,2.$ Then we need to prove
$(\max_{1\leq j \leq p}|z_j|-A_n)/B_n$ converges weakly to the cdf $\exp(-e^{-x})$,
where $A_n=c_n+\frac{1}{2}(1-c_n^2)^{1/2}(n-1)^{-1/2}a_n$,  $B_n=\frac{1}{2}(1-c_n^2)^{1/2}(n-1)^{-1/2}b_n$,
\[
c_n=\Big(\frac{p-1}{n-1}\Big)^{1/2},~~ b_n=b\Big(\frac{nc_n^2}{1-c_n^2}\Big), ~~a_n=a\Big(\frac{nc_n^2}{1-c_n^2}\Big)
\]
with
\begin{equation*}
a(x)=(\log x)^{1/2}-(\log x)^{-1/2}\log(\sqrt{2\pi}\log x)\ \mbox{ and } \ b(x)=(\log x)^{-1/2}
\end{equation*}
for $x>3$. We proceed this  through several steps.

{\it Step 1: Reduction to an easy formulation}. Let $U_i$, $i\ge 1$ be a sequence of i.i.d. random variables uniformly distributed over $(0,1)$, and $U_{1:n}\le U_{2:n}\le \cdots\le U_{n:n}$ be the order statistics of $U_1, U_2, \cdots, U_n$ for each $n\ge 1$. From \eqref{pdf}, the density function of $U_{j: m_n+j-1}$ is
\[
f_{j:m_n+j-1}(x)=\frac{(m_n+j-1)!}{(j-1)!(m_n-1)!}x^{j-1}(1-x)^{m_n-1},~~~~~~x\in (0,1).
\]
Denote the corresponding cdf as $F_{j:m_n+j-1}(x)$. Notice the pdf of $(U_{j: m_n+j-1})^{1/2}$ is proportional to
$x^{2j-1}(1-x^2)^{m_n-1}$. For each $n\geq 2$, let $\{Y_{nj};\, 1\leq j \leq p\}$ be independent random variables such that
$Y_{nj}$ and $(U_{j: m_n+j-1})^{1/2}$ have the same distribution.
By Lemma \ref{independence} and (\ref{elastic}), $\max_{1\leq j \leq p}|z_j|$ and $\max_{1\leq j \leq p}Y_{nj}$ have the same distribution. We claim that, to prove the theorem,  it suffices to show
\begin{eqnarray}\lbl{red_face}
\lim_{n\to\infty}P\Big(\max_{1\leq j \leq p}Y_{nj}^2\le \beta_n(x)\Big)=\exp(-e^{-x})
\end{eqnarray}
for every $x\in \mathbb{R}$, where $\beta_n(x)=c_n^2+c_n(1-c_n^2)^{1/2}(n-1)^{-1/2}(a_n+b_nx)$. In fact, \eqref{red_face}
implies that
\begin{equation}\label{wn}
W_n:=\frac{1}{b_n}\left( c_n^{-1}(1-c_n^2)^{-1/2}(n-1)^{1/2}\big(\max_{1\leq j \leq p}Y_{nj}^2-c_n^2\big)-a_n\right)\td \Lambda,
\end{equation}
where $\Lambda$ is a probability distribution with cdf $\exp(-e^{-x}),\ x\in \mathbb{R}.$
From Taylor's expansion
\begin{eqnarray*}
\max_{1\leq j \leq p}Y_{nj}&=&\left(c_n^2+c_n(1-c_n^2)^{1/2}(n-1)^{-1/2}(a_n+b_nW_n)\right)^{1/2}\\
&=&c_n\left(1+c_n^{-1}(1-c_n^2)^{1/2}\frac{a_n+b_nW_n}{(n-1)^{1/2}}\right)^{1/2} \\
&=&c_n\left(1+\frac12c_n^{-1}(1-c_n^2)^{1/2}\frac{a_n+b_nW_n}{(n-1)^{1/2}}+O_P\Big(\frac{a_n^2}{n-1}\Big)\right)\\
&=&c_n+\frac12(1-c_n^2)^{1/2}(n-1)^{-1/2}a_n+ \frac12(1-c_n^2)^{1/2}(n-1)^{-1/2}b_nW_n\\
& & ~~~~~~~~~~~~~~~~~~~~~~~~~~~~~~~~~~~~~~~~~~~~~~~~~~~~~~~~~~~~~~~~~~~~~+O_P\Big(\frac{\log n}{n}\Big)\\
&=&A_n+B_nW_n+O_p\Big(\frac{\log n}{n}\Big),
\end{eqnarray*}
where we use the facts $a_n \to \infty$, $b_n\to 0$ and $c_n \in (0,1)$ in the above. Since $B_n$ has the scale of $(n\log n)^{-1/2}$, by \eqref{wn},
\[
\frac{\max_{1\leq j \leq p}Y_{nj}-A_n}{B_n}=W_n+O_p\Big(\frac{(\log n)^{3/2}}{n^{1/2}}\Big)\to \Lambda
\]
weakly, which leads to the desired conclusion. Now we proceed to show \eqref{red_face}.

{\it Step 2: A preparation}. We claim that
\begin{equation}\label{monotone2}
1-F_{1:m_n}(x)\le 1-F_{2: m_n+1}(x)\le \cdots\le 1-F_{p: m_n+p-1}(x)
\end{equation}
for $x\in (0,1).$ In fact, since for each $1<j\le p$,
\[
U_{j: m_n+j-1}=\left\{
                 \begin{array}{ll}
                  U_{j-1:m_n+j-2}, & \hbox{ if } U_{m_n+j-1}\le U_{j-1:m_n+j-2};\\
                  \min(U_{j:m_n+j-2}, U_{m_n+j-1}) , & \hbox{ if }  U_{m_n+j-1}>U_{j-1:m_n+j-2},
                 \end{array}
               \right.
\]
which implies that $U_{j-1:m_n+j-2}\le U_{j:m_n+j-1}$  for $1<j\le p$. This yields (\ref{monotone2}).

For each $n\ge 2$, set $a_{nj}=1-F_{p+1-j:m_n+p-j}(\beta_n(x))=1-F_{p+1-j:n-j}(\beta_n(x))$ for $1\le j\le p$.
From \eqref{monotone2}, for each $n$, $a_{ni}$ is non-increasing in $i$.
Since $Y_{nj}^2$ and $U_{j: m_n+j-1}$ are identically distributed, we have
\begin{equation}\label{aprod}
P\Big(\max_{1\leq j \leq p}Y_{nj}^2\le \beta_n(x)\Big)=\prod^p_{j=1}P(Y_{nj}^2\le \beta_n(x))=
\prod^p_{j=1}(1-a_{nj}).
\end{equation}

It is easy to check the following holds: suppose $\{l_n;\, n\geq 1\}$ is sequence of positive integers. Let $z_{ni}\in [0,1)$ be constants for all $1\leq i \leq l_n$ with $\max_{1\le i\le l_n}z_{ni}\to 0$ and
$\sum^{l_n}_{i=1}z_{ni}\to z\in [0,\infty)$.  Then
\begin{equation}\label{lem2}
\lim_{n\to\infty}\prod^{l_n}_{i=1}(1-z_{ni})=e^{-z}.
\end{equation}
Next we will use \eqref{aprod} and \eqref{lem2}  to prove \eqref{red_face}. In fact, we only need to verify that
\begin{equation}\label{sumhasalimita}
\sum^p_{j=1}a_{nj}\to e^{-x}
\end{equation}
and
\begin{equation}\label{maxissmall}
\max_{1\le j\le p}a_{nj}=a_{n1}\to 0.
\end{equation}

{\it Step 3: The analysis of dominated terms}. Fix $\delta\in (\frac12, \frac23)$. Let $j_n=[n^{\delta}]$, the integer part of $n^\delta$. For $1\le j\le j_n$, define
\[
u_{nj}=\frac{(n-j)^{3/2}}{((p-j)m_n)^{1/2}}\left( \beta_n(x)-\frac{p-j}{n-j}\right).
\]
Meanwhile, we rewrite
\[
\beta_n(x)=\frac{p-1}{n-1}+\frac{((p-1)m_n)^{1/2}}{(n-1)^{3/2}}(a_n+b_nx).
\]
Then we see that uniformly over $1\le j\le j_n$,
\beaa
u_{nj}&= & \Big(\frac{p-1}{n-1}-\frac{p-j}{n-j}\Big)\cdot \frac{(n-j)^{3/2}}{((p-j)m_n)^{1/2}}\\
& & ~~~~~~~~~~~~~~~~~~~~~~~~~~~~~~~~~~~~~~ + \Big(\frac{n-j}{n-1}\Big)^{3/2}\cdot \Big(\frac{p-1}{p-j}\Big)^{1/2}(a_n+b_nx)\\
& = & \Big(\frac{p-j}{p-1}\Big)^{-1/2}\cdot \Big(\frac{n-j}{n-1}\Big)^{1/2}\cdot
\Big(\frac{n-p}{p-1}\Big)^{1/2}\cdot\Big(\frac{n-1}{n}\Big)^{-1/2}\cdot \frac{j-1}{n^{1/2}}\\
& & ~~~~~~~~~~~~~~~~~~~~~~~~~~~~~~~~~~~~~~ + \Big(\frac{n-j}{n-1}\Big)^{3/2}\cdot \Big(\frac{p-j}{p-1}\Big)^{-1/2}(a_n+b_nx).
\eeaa
Now, $\frac{n-p}{p-1}=\frac{1-c_n^2}{c_n^2}$. Also, given $\tau\in \mathbb{R}$, trivially $\Big(\frac{p-j}{p-1}\Big)^{\tau}=1+O(\frac{j}{n})$ and $\Big(\frac{n-j}{n-1}\Big)^{\tau}=1+O(\frac{j}{n})$ uniformly for all $1\leq j \leq j_n.$ Since $a_n \sim (\log n)^{1/2}$ and $b_n=o(1)$, we have
\begin{eqnarray}\label{onlyone}
u_{nj} &=& \frac{(1-c_n^2)^{1/2}}{n^{1/2}c_n}(j-1)(1+o(1)) + a_n+ b_nx + O\Big(\frac{j\log n}{n}\Big)\nonumber\\
& = & \frac{(1-c_n^2)^{1/2}}{n^{1/2}c_n}(j-1)(1+o(1)) + a_n+ b_nx
\end{eqnarray}
uniformly for all $1\leq j \leq j_n.$

In Lemma~\ref{reiss}, take $r=n-j$ and $k=n-p$ to have
\[
\sup_{B\in \cal{B}}\Big|P(V_{p-j+1:n-j}\in B)-\int_B(1+l_1(t)+l_2(t))\phi(t)dt\Big|=O(n^{-3/2})
\]
uniformly over $1\le j\le j_n$ as $n\to\infty$,  where
\[
V_{p-j+1:n-j}=\frac{(n-j)^{3/2}}{((n-j)(n-p))^{1/2}}\Big(U_{p-j+1:n-j}-\frac{p-j}{n-j}\Big)
\]
and
where,
for $i=1,2$, $l_i(t)$ is a polynomial in $t$ of degree $\le 3i$, depending on $n$,  and all of its coefficients  are of order $O(n^{-i/2})$ by the assumption $h_1<\frac{p}{n} < h_2$ for all $n\geq 2.$ Now, by taking $B=(u_{nj}, \infty)$ we obtain
\[
a_{nj}=P(V_{p-j+1:n-j}>u_{nj})
=\int^\infty_{u_{nj}}(1+l_1(t)+l_2(t))\phi(t)dt+O(n^{-3/2})
\]
uniformly for $1\le j\le j_n$ as $n\to\infty$. From L'Hospital's rule, we have  that for any $r\ge 0$
\begin{equation}\label{tailapp}
\int^\infty_xt^r\phi(t)dt\sim x^{r-1}\phi(x)~~~~~\mbox{ as } x\to\infty.
\end{equation}
Since $\min_{1\le j\le j_n}u_{nj}\to\infty$ as $n\to\infty$ by (\ref{onlyone}), it  follows from \eqref{tailapp} that
\[
\int^\infty_{u_{nj}}t^r\phi(t)dt\sim (u_{nj})^{r-1}\phi(u_{nj}) =
\left\{
  \begin{array}{ll}
   \frac{\phi(u_{nj})}{u_{nj}} , & \hbox{ if } r=0;\\
   O((\max_{1\le j\le j_n}u_{nj})^r)\frac{\phi(u_{nj})}{u_{nj}} , & \hbox{ if } r>0
  \end{array}
\right.
\]
holds uniformly over $1\le j\le j_n$. Furthermore, since the coefficients of
$l_i(t)$ are uniformly bounded by $O(n^{-i/2})$ for $i=1,2$,
we have
\begin{eqnarray*}
&&\int^\infty_{u_{nj}}(1+l_1(t)+l_2(t))\phi(t)dt\\
&=&\Big[1+O\Big(\frac{(\max_{1\le j\le j_n}u_{nj})^3}{n^{1/2}}\Big)+O\Big(\frac{(\max_{1\le j\le j_n}u_{nj})^6}{n}\Big)\Big]\frac{\phi(u_{nj})}{u_{nj}}\\
&=&(1+O(\frac{j_n^3}{n^2}))\frac{\phi(u_{nj})}{u_{nj}}
\end{eqnarray*}
uniformly  over $1\le j\le j_n$, and thus obtain that
\begin{equation}\label{pestimate}
a_{nj}=\left(1+O(\frac{j_n^3}{n^2})\right)\frac{\phi(u_{nj})}{u_{nj}}+O(n^{-3/2})
\end{equation}
uniformly  over $1\le j\le j_n$.  Therefore, we have
\begin{equation}\label{sum1tojn}
\sum^{j_n}_{j=1}a_{nj}
=(1+o(1))\sum^{j_n}_{j=1}\frac{\phi(u_{nj})}{u_{nj}}+o(1).
\end{equation}

In Lemma~\ref{seed}, by taking $x_n=nc_n^2/(1-c_n^2)$ and $c_{nj}=x_n^{-1/2}(1+o(1))$ where ``$o(1)$" is as indicated in (\ref{onlyone}), we then get
\begin{equation}\label{sum1}
\lim_{n\to\infty}\sum^{j_n}_{j=1}\frac{\phi(u_{nj})}{u_{nj}}= e^{-x}.
\end{equation}

{\it Step 4: Non-dominated terms are negligible}.
From \eqref{onlyone} again, we see
\[
u_{nj_n}^2\ge 2c_n^{-1}(1-c_n^2)^{1/2}(a_n+b_nx)\frac{n^{\delta}-2}{n^{1/2}}(1+o(1))\ge 6\log n
\]
for all large $n$. Hence,
\[
\phi(u_{nj_n})\le \frac{1}{\sqrt{2\pi}}\exp(-3\log n)=\frac{1}{\sqrt{2\pi}n^{3}}
\]
for all large $n$. Then it follows from \eqref{pestimate} that $a_{nj_n}=O(n^{-3/2})$, and hence
\[
\sum^p_{j=j_n+1}a_{nj}\le (n-j_n)a_{nj_n}=O(n^{-1/2}).
\]
This together with \eqref{sum1}  yields \eqref{sumhasalimita}.
The proof is then completed. \hfill $\blacksquare$

\subsection{The Proof of Theorem \ref{Theorem_product}}\lbl{jinren}

We begin with some preparation.  The following result characterizes the structure of the radius of the eigenvalues from the product ensemble.
\begin{lemma}\lbl{Catroll} Let $k$ and $z_1,\cdots, z_n$ be as in (\ref{up_stairs}). Let $\{s_{j, r},\,1\le r\le k,  j\ge 1\}$ be independent random variables and $s_{j, r}$ have the Gamma  density $y^{j-1}e^{-y}I(y>0)/(j-1)!$ for each $j$ and $r.$
Then $\max_{1\leq j \leq n}|z_j|^2$ and $\max_{1\leq j \leq n}\prod_{r=1}^ks_{j,r}$ have the same distribution.
\end {lemma}
\textbf{Proof}.
Let $\{s_{j, r};\, 1\le r\le k,  j\ge 1\}$ be independent random variables and $s_{j, r}$
follow a Gamma($j$) distribution with density function $y^{j-1}e^{-y}I(y\geq 0)/\Gamma(j)$ for all $1\le r\le k$ and  $j\ge 1.$ Define $v_1(y)=\exp(-y)$, $y>0$, and set for $j\ge 2$
\bea\lbl{blackmail}
v_j(y)=\int^\infty_{0}v_{j-1}(y/s)\frac{e^{-s}}{s}ds.
\eea
One can easily verify that for each $j\ge 1$, $v_j(y)$ is proportional to $w_j(y^{1/2})$, i.e.,
for some constants $d_j>0$,
\begin{equation}\label{fact-1}
w_j(y^{1/2})=d_jv_j(y), ~~~~y>0.
\end{equation}
Let $z$ be any complex number with $\textit{Re}(z)>0$,  and define for $j\ge 1$
\[
\gamma_j(z)=\int^\infty_0y^{z-1}v_j(y)dy.
\]
Note that $\gamma_1(z)=\Gamma(z)=\int_0^\infty y^{z-1}e^{-y}dy$. For $j\ge 2$, by using (\ref{blackmail}),
\beaa
\gamma_j(z) & = & \int_0^{\infty}\frac{e^{-s}}{s}\Big[\int_0^{\infty}y^{z-1}v_{j-1}\big(\frac{y}{s}\big)\,dy\Big]\,ds\\
& = & \int_0^{\infty}s^{z-1}e^{-s}\Big[\int_0^{\infty}y^{z-1}v_{j-1}(y)\,dy\Big]\,ds=\Gamma(z)\gamma_{j-1}(z).
\eeaa
 Thus, we have
\begin{equation}\label{fact-2}
\gamma_{j}(z)=\Gamma(z)^{j},~~~~j\ge 1.
\end{equation}

Assume $Y_{nj}$,  $1\le j\le n$ are independent random variables,  and for each $1\le j\le n$, the density of $Y_{nj}$ is  proportional to $y^{2j-1}w_k(y)$.  By Lemma \ref{independence} and (\ref{elastic}), $\max_{1\le j\le n}|z_j|^2$ and
$\max_{1\le j\le n}Y_{nj}^2$ are identically distributed. Furthermore, since the density function of $Y_{nj}^2$, denoted by $f_j(y)$,  is proportional to $y^{j-1}w_k(y^{1/2})$, and thus proportional to
 $y^{j-1}v_k(y)$ from \eqref{fact-1},  we have from \eqref{fact-2} that
\[
f_j(y)=\frac{y^{j-1}v_k(y)}{\int^\infty_0y^{j-1}v_k(y)dy}
=\frac{y^{j-1}v_k(y)}{\Gamma(j)^k},  ~~~y>0,
\]
for $1\le j\le n$. Let the characteristic function of $\log Y_{nj}^2$ be denoted by $g_j(t)$. Then
we have
\[
g_j(t)=\frac{1}{\Gamma(j)^k}\int^\infty_0e^{it\log y}y^{j-1}v_k(y)dy
=\frac{1}{\Gamma(j)^k}\int^\infty_0y^{j-1+it}v_k(y)dy=
\left(\frac{\Gamma(j+it)}{\Gamma(j)}\right)^k
\]
from  \eqref{fact-2}. Since $\Gamma(j+it)/\Gamma(j)$ is the characteristic function
of $\log s_{j, r}$, it follows that  $\log Y_{nj}^2$ has the same distribution as that of
$\sum^{k}_{r=1}\log s_{j, r}$, or equivalently, $Y_{nj}^2$ has the same distribution as that of
$\prod^{k}_{r=1}s_{j, r}$ for $j\ge 1$. This implies the desired conclusion.      \hfill$\blacksquare$


\begin{lemma}\lbl{max-in} Let $k$  be as in (\ref{up_stairs}) and $\{s_{j, r},\,1\le r\le k,  j\ge 1\}$ be independent r.v.'s such that $s_{j, r}$ has density $y^{j-1}e^{-y}I(y>0)/(j-1)!$ for all $j,r$. Set $\eta(x)=x-1-\log x$ and
\bea\lbl{dinner_mala}
M_n(i)=\max_{n-i+1\le j\le n}\Big|\sum^{k}_{r=1}\Big(\eta\Big(\frac{s_{j,r}}{j}\Big)-E\eta\Big(\frac{s_{j,r}}{j}\Big)\Big)\Big|, \ \ 1\leq i \leq n.
\eea
Set $\psi(x)=\frac{\Gamma'(x)}{\Gamma(x)}$ for $x>0.$ Then for $1\leq i \leq n$
\beaa
\Big|\max_{n-i+1\le j\le n}\log \prod^{k}_{r=1} s_{j, r}-\max_{n-i+1\le j\le n}\Big(\frac{1}{j}\sum^{k}_{r=1}(s_{j,r}-j)+k\,\psi(j)\Big)\Big|
\le M_n(i).
\eeaa
\end{lemma}
\textbf{Proof}. Set $Y_j=\prod^{k}_{r=1}s_{j, r}$ for $j\geq 1.$ Then,
\[
\log Y_j=\sum^{k}_{r=1}\log s_{j, r}
\]
for $j\ge 1$. The moment generating functions of $\log s_{j,r}$ is
\begin{equation}\label{ymoment}
m_j(t)=E(e^{t\log s_{j,r}})=\frac{\Gamma(j+t)}{\Gamma(j)}
\end{equation}
for $t>-j.$ Therefore,
\[
E(\log s_{j,r})=\frac{d}{dt}m_j(t)|_{t=0}=\frac{\Gamma'(j)}{\Gamma(j)}=\psi(j)
\]
by (\ref{maomi}). Note that $\eta(x)=x-1-\log x$ for $x>0$. Since $\eta(x)=\int^x_1\frac{s-1}{s}ds$, it is easy to verify that
\begin{equation}\label{eta-est}
0\le \eta(x)\le \frac{(x-1)^2}{2\min(x,1)}, ~~~~~~x>0.
\end{equation}
By using the expression $\log x=x-1-\eta(x)$ we can rewrite $\log Y_{j}$ as
\begin{eqnarray*}
\log Y_j&=&\sum^{k}_{r=1}\log \frac{s_{j,r}}{j}+k\log j\\
&=&\sum^{k}_{r=1}\frac{s_{j,r}-j}{j}-\sum^{k}_{r=1}\eta\Big(\frac{s_{j,r}}{j}\Big)+k\log j\\
&=&\frac{1}{j}\sum^{k}_{r=1}(s_{j,r}-j)-\sum^{k}_{r=1}\eta\Big(\frac{s_{j,r}}{j}\Big)+k\log j\\
&=&\frac{1}{j}\sum^{k}_{r=1}(s_{j,r}-j)+k\psi(j)-\sum^{k}_{r=1}\eta\Big(\frac{s_{j,r}}{j}\Big)+k(\log j-\psi(j)).
\end{eqnarray*}
Since $E(\log Y_j)=k\psi(j)$,  we see that
\[
\sum^{k}_{r=1}E\eta\Big(\frac{s_{j,r}}{j}\Big)=k(\log j-\psi(j))
\]
and thus
\begin{equation}\label{logY}
\log Y_j
=\frac{1}{j}\sum^{k}_{r=1}(s_{j,r}-j)+k\psi(j)-\sum^{k}_{r=1}\Big(\eta\Big(\frac{s_{j,r}}{j}\Big)-E\eta\Big(\frac{s_{j,r}}{j}\Big)\Big).
\end{equation}
Note that for any two sequences of reals numbers $\{x_n\}$ and $\{y_n\}$,
\beaa
|\max_{1\le j\le n}x_j-\max_{1\le j\le n}y_j|\le \max_{1\le j\le n}|x_j-y_j|.
\eeaa
Then it follows from \eqref{logY} that
\beaa
\ \ \ \ \ \ \ \Big|\max_{n-i+1\le j\le n}\log Y_j-\max_{n-i+1\le j\le n}\Big(\frac{1}{j}\sum^{k}_{r=1}(s_{j,r}-j)+k\psi(j)\Big)\Big|
\le M_n(i). ~~~~~~~~~~~ \blacksquare\\
\eeaa

We estimate $M_n(\cdot)$ next.
\begin{lemma}\label{smallMM} Let $k$  be as in (\ref{up_stairs}) and $M_n(i)$ be defined as in Lemma \ref{max-in}. Assume $\{j_n;\, n\geq 1\}$ is a sequence of numbers satisfying $1\le j_n\le \frac{1}{2}n$ for all $n$. Then, for any sequence of positive integers $\{k_n\}$, $M_n(j_n)=O_P\Big(\frac{j_nk_n^{1/2}}{n}\Big).$ Further, if $\lim_{n\to\infty}k_n/n= 0$, then $M_n(j_n)=O_P\Big(\frac{k_n\log n}{n}\Big).$
\end{lemma}

\noindent\textbf{Proof.}
By using the Minkowski inequality and \eqref{eta-est} we get
\begin{eqnarray*}
E(M_n(j_n))&\leq &\sum_{n-j_n+1\le j\le n}E\sum^{k_n}_{r=1}\Big|\eta\Big(\frac{s_{j,r}}{j}\Big)-E\eta\Big(\frac{s_{j,r}}{j}\Big)\Big|\\
&\le&\sum_{n-j_n+1\le j\le n}\Big[E\Big(\sum^{k_n}_{r=1}\Big|\eta\Big(\frac{s_{j,r}}{j}\Big)-E\eta\Big(\frac{s_{j,r}}{j}\Big)\Big|\Big)^2\Big]^{1/2}\\
&\leq &\sum_{n-j_n+1\le j\le n}\Big[\sum^{k_n}_{r=1} E\Big|\eta\Big(\frac{s_{j,r}}{j}\Big)-E\eta\Big(\frac{s_{j,r}}{j}\Big)\Big|^2\Big]^{1/2}\\
&\le&k_n^{1/2}\sum_{n-j_n+1\le j\le n}\Big(E\eta\Big(\frac{s_{j,1}}{j}\Big)^2\Big)^{1/2}\\
&\le& \frac12 k_n^{1/2}\sum_{n-j_n+1\le j\le n}\Big\{
E\Big[\Big(\frac{s_{j,1}-j}{j}\Big)^4\Big(\min\Big(\frac{s_{j,1}}{j}, 1\Big)\Big)^{-2}\Big]\Big\}^{1/2}
\end{eqnarray*}
by (\ref{eta-est}). Since $s_{j,1}$ has  density $y^{j-1}e^{-y}I(y>0)/(j-1)!$, we see that  $E\big(s_{j,1}^{-4}\big)=\frac{\Gamma(j-4)}{\Gamma(j)}$. By the Marcinkiewicz-Zygmund inequality (see, for example, Corollary 2 from Chow and Teicher, 2003), we obtain $E(s_{j,1}-j)^8\le Kj^4$ for any $j\ge 1$ where $K$ is a constant not depending on $j$. Then, it follows from H\"older's inequality that
\begin{eqnarray*}
& & E\Big(\Big(\frac{s_{j,1}-j}{j}\Big)^4\Big(\min\Big(\frac{s_{j,1}}{j}, 1)\Big)^{-2}\Big)\\
&\le&
\Big[E\Big(\frac{s_{j,1}-j}{j}\Big)^8\cdot E\Big(\min\Big(\frac{s_{j,1}}{j}, 1\Big)\Big)^{-4}\Big]^{1/2}\\
&\le&\left[E\Big(\frac{s_{j,1}-j}{j}\Big)^8\cdot E\Big(\Big(\frac{j}{s_{j,1}}\Big)^4+1\Big)\right]^{1/2}\\
&\le&\Big[\Big(\frac{j^3}{(j-1)(j-2)(j-3)}+1\Big)^{1/2}\Big(E\Big(\frac{s_{j,1}-j}{j}\Big)^8\Big]^{1/2}\\
&\le& Cj^{-2}
\end{eqnarray*}
 for any $j\geq 4$ where $C$ is a constant.
Combining the last two assertions, we get $E(M_n(j_n)) \leq O(\frac{j_nk_n^{1/2}}{n}).$ This implies the first conclusion.


Now we prove the second one. Recall $\psi(x)=\frac{\Gamma'(x)}{\Gamma(x)}$ for $x>0$ as in (\ref{maomi}). By Formulas 6.3.18 and 6.4.12
from Abramowitz and Stegun (1972),
\begin{equation}\label{psi}
\psi(x)=\log x-\frac{1}{2x}+O\Big(\frac{1}{x^2}\Big)\ \ \mbox { and }\ \ \psi'(x)=\frac{1}{x}+\frac{1}{2x^2}+O\Big(\frac{1}{x^3}\Big)
\end{equation}
as $x\to +\infty$. It is easy to check $E\log s_{j,1}= \frac{1}{\Gamma(j)}\int_0^{\infty}(\log y)y^{j-1}e^{-y}\,dy=\psi(j)$. Thus, from the first expression, we have
\[
E\eta\Big(\frac{s_{j,1}}{j}\Big)=\log j-\psi(j)=O\Big(\frac{1}{j}\Big)
\]
as $j\to\infty$. Hence, by (\ref{dinner_mala}),
\begin{equation}\label{M-estimation}
M_n(j_n)\le\max_{n-j_n+1\le j\le n}\sum^{k_n}_{r=1}\eta\Big(\frac{s_{j,r}}{j}\Big)+ O\Big(\frac{k_n}{n}\Big).
\end{equation}
By Theorem 1 on page 217 from  Petrov (1975), we have that
\begin{equation}\label{normapp12}
P(s_{j,1}>j+ j^{1/2}x)=(1+o(1))(1-\Phi(x))
\end{equation}
uniformly for $x\in (0, a_n)$ and  $n/2\le j\le n$ as $n\to\infty$, where $\{a_n;\, n\geq 1\}$ is an arbitrarily given sequence of positive numbers with $a_n=o(n^{1/6})$. By taking $r=0$ in (\ref{tailapp}), we see that $1-\Phi(x) \sim \frac{1}{\sqrt{2\pi}\,x}e^{-x^2/2}$ as $x\to +\infty.$ Now select $x=2(\log n)^{1/2}$ in (\ref{normapp12}) to have
\[
P\big(s_{j,1}>j+ 2j^{1/2}(\log n)^{1/2}\big)=(1+o(1))\big(1-\Phi(2(\log n)^{1/2}\big)=O(\frac{1}{n^2})
\]
uniformly for $n/2\le j\le n$ as $n\to\infty$. Similarly we have
\[
P\big(s_{j,1}<j-2j^{1/2}(\log n)^{1/2}\big)=(1+o(1))\big(1-\Phi(2(\log n)^{1/2}\big)=O\Big(\frac{1}{n^2}\Big)
\]
uniformly for $n/2\le j\le n$ as $n\to\infty$. This implies
\[
\sum^{k_n}_{r=1}\sum^n_{j=n-j_n+1}P\big(|s_{j,r}-j|>2j^{1/2}(\log n)^{1/2}\big)=O\Big(\frac{j_nk_n}{n^2}\Big)=O\Big(\frac{k_n}{n}\Big)=o(1),
\]
and thus we get
\[
\max_{1\le r\le k_n}\max_{n-j_n+1\le j\le n}\Big|\frac{s_{j,r}}{j}-1\Big|=O_P\Big(\frac{(\log n)^{1/2}}{n^{1/2}}\Big).
\]
Consequently,
\[
\min_{1\le r\le k_n}\min_{n-j_n+1\le j\le n}\frac{s_{j,r}}{j}=1+O_P\Big(\frac{(\log n)^{1/2}}{n^{1/2}}\Big).
\]
By \eqref{M-estimation} and then \eqref{eta-est}, we obtain
\begin{eqnarray*}
M_n(j_n)&\le& k_n\cdot\max_{n-j_n+1\le j\le n}\max_{1\le r\le k_n}\eta(\frac{s_{j,r}}{j})+ O(\frac{k_n}{n})\\
&\le& \frac{k_n}2\cdot\frac{\max_{n-j_n+1\le j\le n}\max_{1\le r\le k_n}
|\frac{s_{j,r}}{j}-1|^2}{\min\{1,\min_{n-j_n+1\le j\le n}\min_{1\le r\le k_n}\frac{s_{j,r}}{j}\}}+O\Big(\frac{k_n}{n}\Big)\\
&=&O_P\Big(\frac{k_n\log n}{n}\Big),
\end{eqnarray*}
proving the second conclusion.  \hfill$\blacksquare$\\

Review the notation we use before: $\psi(x)=\frac{\Gamma'(x)}{\Gamma(x)}$ for $x>0$ as in (\ref{maomi}) and
\bea\label{no_xian}
Y_j=\prod^{k}_{r=1}s_{j, r}
\eea
for $j\geq 1$, where
$\{s_{j, r},\,1\le r\le k,  j\ge 1\}$ are independent random variables such that $s_{j, r}$ has density $y^{j-1}e^{-y}I(y>0)/(j-1)!$ for all $j,r$.
\begin{lemma}\label{probability4smallj}
Let $\{j_n;\, n\geq 1\}$ and $\{k_n;\, n\geq 1\}$ be  positive integers satisfying $\lim_{n\to\infty}\frac{j_n}n=0$ and
$\lim_{n\to\infty}(\frac{k_n}{n})^{1/2}\frac{j_n}{(\log n)^{1/2}}=\infty.$  Then for any $x\in \mathbb{R}$,
\begin{equation}\label{sumissmall}
\lim_{n\to\infty}\sum^{n-j_n}_{j=1}P\Big(\log Y_j>k_n\psi(n)+ \Big(\frac{k_n}{n}\Big)^{1/2}x\Big)=0.
\end{equation}
\end{lemma}

\noindent\textbf{Proof.} Fix $x\in \mathbb{R}$. It follows from \eqref{ymoment} that for each $1 \le j\le n-j_n$ and any $t>0$,
\begin{eqnarray*}
&&P\Big(\log Y_j>k_n\psi(n)+\Big(\frac{k_n}{n}\Big)^{1/2}x\Big)\\
&\le& \frac{E(e^{t\log Y_j})}{\exp\big\{t(k_n\psi(n)+(\frac{k_n}{n})^{1/2} x)\big\}}\\
&=& \exp\Big\{k_n(\log(\Gamma(j+t)-\log\Gamma(j)))-t\Big(k_n\psi(n)+\Big(\frac{k_n}{n}\Big)^{1/2}x\Big)\Big\}\\
&=& \exp\Big\{k_n\int^t_0\psi(j+s)ds-t\Big(k_n\psi(n)+\Big(\frac{k_n}{n}\Big)^{1/2}x\Big)\Big\}\\
&=& \exp\Big\{k_n\int^t_0[\psi(j+s)-\psi(j)]ds-t\Big[k_n(\psi(n)-\psi(j))+\Big(\frac{k_n}{n}\Big)^{1/2}x\Big]\Big\}.
\end{eqnarray*}
From \eqref{psi}, there exist an integer $j_0$ such that for all $j_0\le j\le n-j_n$
\[
\log\frac{j+s}{j}\le \psi(j+s)-\psi(j)=\int^{s}_0\psi'(j+v)dv\le \frac{1.1s}{j}, ~~~s\ge 0.
\]
By the first inequality above, for
all large $n$,
\[
\psi(n)-\psi(j)\ge \log\frac{n}{j} \ge \log\frac{n}{n-j_n}=-\log(1-\frac{j_n}{n})\ge \frac{0. 999j_n}{n},  ~~~j_0\le j\le n-j_n.
\]
Hence, by assumption $(\frac{k_n}{n})^{1/2}=o(\frac{j_nk_n}{n})$  we see that
\[
k_n(\psi(n)-\psi(j))+\Big(\frac{k_n}{n}\Big)^{1/2}x\ge 0.99k_n\log \frac{n}{j},~~~~j_0\le j\le n-j_n
\]
for all large $n$. Therefore we have for $j_0\le j\le n-j_n$
\begin{eqnarray*}
&&P\Big(\log Y_j>k_n\psi(n)+\Big(\frac{k_n}{n}\Big)^{1/2}x\Big)\\
&\le& \exp\Big\{1.1k_n\int^t_0\frac{s}{j}ds-0.99tk_n(\log n-\log j)\Big\}\\
&=& \exp\Big\{k_n\Big(\frac{1.1t^2}{2j}-0.99t(\log n-\log j)\Big)\Big\}
\end{eqnarray*}
for all $t>0$ and large $n$ which does not depend on $t.$ By selecting $t=0.99j(\log n-\log j)$ we have
\bea
& & P\Big(\log Y_j>k_n\psi(n)+\Big(\frac{k_n}{n}\Big)^{1/2}x\Big) \nonumber\\
& \le & \exp\big\{-0.44k_nj(\log n-\log j)^2\big\},~~~~j_0\le j\le n-j_n, \lbl{Yuan}
\eea
for all large $n$.   Note that
\begin{eqnarray*}
\min_{j_0\le j\le n-j_n}j(\log n-\log j)^2
&\ge& \min_{j_0\le s\le n-j_n}s(\log n-\log s)^2\\
&=&\min_{j_0^{1/2}\le s\le (n-j_n)^{1/2}}s^2(\log n-2\log s)^2\\
&=&\Big(\min_{j_0^{1/2}\le s\le (n-j_n)^{1/2}}s(\log n-2\log s)\Big)^2,
\end{eqnarray*}
where the last three minima are taken over all real numbers satisfying the corresponding constraints. It is easily seen that the minimum of $s(\log n-2\log s)$ for $j_0^{1/2}\le s\le (n-j_n)^{1/2}$ is achieved at the two end points
of the interval, $t=j_0^{1/2}$ or $s=(n-j_n)^{1/2}$. Thus, for all large $n$,
 \begin{eqnarray*}
\min_{j_0\le j\le n-j_n}j(\log n-\log j)^2
&\ge & \min\big\{j_0(\log n-\log j_0)^2, (n-j_n)(\log n-\log (n-j_n))^2\big\} \\
&\ge &\min\Big\{\frac12(\log n)^2,\, \frac12\frac{j_n^2}{n}\Big\}.
\end{eqnarray*}
From the given condition $(\frac{k_n}{n})^{1/2}\frac{j_n}{(\log n)^{1/2}}=\infty$, we obtain
\[
k_n\min_{j_0\le j\le n-j_n}j(\log n-\log j)^2\ge 10\log n
\]
for all large $n$. Therefore, combining all of the inequalities from (\ref{Yuan}) to the above, we have
\[
\max_{j_0\le j\le n-j_n}P\Big(\log Y_j>k_n\psi(n)+\Big(\frac{k_n}{n}\Big)^{1/2}x\Big)\le \exp(-4.4\log n)=n^{-4.4},
\]
and hence
\[
\sum^{n-j_n}_{j=j_0}P\Big(\log Y_j>k_n\psi(n)+\Big(\frac{k_n}{n}\Big)^{1/2}x\Big)=O(n^{-3.4})\to 0.
\]
Finally, observe that, for each $1\le j <j_0$, $\log Y_j$ is a sum of $k_n$'s many i.i.d. random variables with $Ee^{t\log Y_j}<\infty$ for all $|t|<\frac{1}{2}$. Then, by the Chernoff bound (see, for instance, p. 27 from Dembo and Zeitouni, 1998),
\[
P\Big(\log Y_j>k_n\psi(n)+\Big(\frac{k_n}{n}\Big)^{1/2}x\Big)\to 0.
\]
The last two assertions imply the desired result. \hfill$\blacksquare$
\vspace{20pt}


Recall $\Lambda(x)=\exp(-e^{-x})$ for all $x\in \mathbb{R}.$ Considering convenience of formulation, we first prove the following proposition from which
Theorem \ref{Theorem_product} will be obtained.

\begin{prop}\label{prop5} Let $\psi(x)$ be as in (\ref{maomi}), $a(x)$ and $b(x)$  be as in Theorem \ref{prop5}, and $z_j$'s and $k_n$ be as in Theorem \ref{Theorem_product}. Define $\Phi_0(y)=\Lambda(y)$,  $a_n=a(n/k_n)$, $b_n=b(n/k_n)$ if $\alpha=0$, and $a_n=0$,  $b_n=1$
if $\alpha\in (0, \infty]$. Then
\begin{equation}\label{step0}
\lim_{n\to\infty}P\big(\frac{\max_{1\le j\le n}\log|z_j|-k_n\psi(n)/2}{(k_n/n)^{1/2}/2}\le
a_n+b_ny\big)=\Phi_{\alpha}(y),~~~y\in \mathbb{R}.
\end{equation}
\end{prop}


\noindent{\bf Proof.}
For each of the three cases: $\alpha=0$, $\alpha\in (0,\infty)$, and $\alpha=\infty$ we will show that
there exists a sequence of positive integers $\{j_n\}$ with $1\le j_n\le n/2$ such that
\bea
&&\lim_{n\to\infty}\sum^{n-j_n}_{j=1}P\Big(\log Y_j>k_n\psi(n)+\Big(\frac{k_n}{n}\Big)^{1/2}(a_n+b_ny)\Big)=0,
 ~~~  y\in \mathbb{R};\label{stepone}\\
&&\frac{M_n(j_n)}{(k_n/n)^{1/2}b_n}~~\text{converges in probability to zero}\label{steptwo}
\eea
where $M_n(\cdot)$ is defined as in Lemma \ref{max-in}, and
\bea\label{stepthree}
& & \lim_{n\to\infty}P\Big(\frac{\max\limits_{n-j_n+1\le j\le n}
\Big(\frac{1}{j}\sum^{k_n}_{r=1}(s_{j,r}-j)
+k_n\psi(j)\Big)-k_n\psi(n)}
{(k_n/n)^{1/2}}\le a_n+b_ny\Big)  ~~~~~~~~~~~~~~~~ \nonumber\\
&=& \Phi_\alpha(y)
\eea
for $y\in \mathbb{R}$, where
$\{s_{j, r},\,1\le r\le k,  j\ge 1\}$ are independent random variables such that $s_{j, r}$ has density $y^{j-1}e^{-y}I(y>0)/(j-1)!$ for all $j$ and $r$.  In fact, \eqref{stepthree} implies
\[
\frac{\max\limits_{n-j_n+1\le j\le n}\big(\frac{1}{j}\sum^{k_n}_{r=1}(s_{j,r}-j)+k_n\psi(j)\big)-k_n\psi(n)}
{(k_n/n)^{1/2}b_n}-\frac{a_n}{b_n}\td \Phi_\alpha.
\]
Review the definition of $Y_j$ in (\ref{no_xian}). The above result together with \eqref{steptwo}, Lemmas \ref{Catroll} and \ref{max-in} implies that
\[
\frac{\max\limits_{n-j_n+1\le j\le n}\log Y_j-k_n\psi(n)}
{(k_n/n)^{1/2}b_n}-\frac{a_n}{b_n}\td \Phi_\alpha.
\]
Since \eqref{stepone} implies
\bea
\frac{\max\limits_{1\le j\le n-j_n}\log Y_j-k_n\psi(n)}
{(k_n/n)^{1/2}b_n}-\frac{a_n}{b_n}\td 0,
\eea
the two limits above imply \eqref{step0} due to the fact that $\max_{1\le j\le n}\log |z_j|$ and $\frac12\max_{1\le j\le n}\log Y_j$ are
identically distributed by Lemma \ref{Catroll}.



Now we start to verify equations \eqref{stepone}-\eqref{stepthree} with a choice of $j_n$ given by
\bea\label{onjn}
j_n=\mbox{ the integer part of } \Big(\frac{n}{k_n}\Big)^{1/2}n^{1/8}+1
\eea
for all large $n$.


 {\it Proof of \eqref{stepone}}. It is easy to verify that the conditions in Lemma~\ref{probability4smallj} are satisfied, and thus \eqref{sumissmall} holds. In case $\alpha\in (0, \infty]$,  $a_n=0$ and  $b_n=1$, and \eqref{stepone} holds in this case. When $\alpha=0$, $a_n+b_ny>0$ for all large $n$, by applying \eqref{sumissmall} with $x=0$ we have
\beaa
& & \lim_{n\to\infty}\sum^{n-j_n}_{j=1}P\Big(\log Y_j>k_n\psi(n)+\Big(\frac{k_n}{n}\Big)^{1/2}(a_n+b_ny)\Big)\\
&\le &
\lim_{n\to\infty}\sum^{n-j_n}_{j=1}P(\log Y_j>k_n\psi(n))=0,
\eeaa
that is, \eqref{stepone} holds.  This completes the proof of \eqref{stepone} for all three cases.

 {\it Proof of \eqref{steptwo}}. To prove \eqref{steptwo}, it suffices to show $M_n(j_n)=O_P((\frac{k_n}{n})^{1/2}(\log n)^{-1})$ since
since $b_n \geq (\log n)^{-1/2}$ for all large $n.$ We use
Lemma~\ref{smallMM} this time. When $\alpha\in (0, \infty]$, $j_n=O_P(n^{1/8})$ from \eqref{onjn},
and then we have from the first conclusion in Lemma \ref{smallMM}  that
\[
M_n(j_n)=O_P\Big(\frac{k_n^{1/2}}{n}j_n\Big)=O_P\Big(\Big(\frac{k_n}{n}\Big)^{1/2}n^{-3/8}\Big)=O_P\Big(\Big(\frac{k_n}{n}\Big)^{1/2}(\log n)^{-1}\Big).
\]
When $\alpha=0$, we have from the two conclusions in Lemma \ref{smallMM}  that
\beaa
M_n(j_n) & = & O_P\Big(\min\Big\{\frac{j_nk_n^{1/2}}{n}, \frac{k_n\log n}{n}\Big\}\Big)\\
&=& \Big(\frac{k_n}{n}\Big)^{1/2}\cdot O_P\Big(\min\Big\{\frac{n^{1/8}}{k_n^{1/2}}, \Big(\frac{k_n}{n}\Big)^{1/2}\log n\Big\}\Big)\\
& = & \Big(\frac{k_n}{n}\Big)^{1/2}\cdot O_P\big(n^{-1/8}\big)\\
&=& O_P\Big(\Big(\frac{k_n}{n}\Big)^{1/2}(\log n)^{-1}\Big)
\eeaa
since $\frac{n^{1/8}}{k_n^{1/2}} \le n^{-1/8}$ if $k_n\ge n^{1/2}$ and $\frac{k_n^{1/2}\log n}{n^{1/2}}\le n^{-1/8}$
 if $k_n< n^{1/2}$.

 {\it Proof of \eqref{stepthree}}. Set $T_n(j_n)=\max\limits_{n-j_n+1\le j\le n}\big(\frac{1}{j}\sum^{k_n}_{r=1}
 (s_{j,r}-j)+k_n\psi(j)\big)$. Then
\bea
&&P\Big(T_n(j_n)\le k_n\psi(n)+\Big(\frac{k_n}{n}\Big)^{1/2}(a_n+b_ny)\Big) \lbl{noxing}\\
 &=& \prod^{n}_{j=n-j_n+1}P\Big(\sum^{k_n}_{r=1}
 (s_{j,r}-j) \leq jk_n(\psi(n)-\psi(j))+ j\Big(\frac{k_n}{n}\Big)^{1/2}(a_n+b_ny)\Big). \lbl{RiemannGeo}
\eea
Notice $\sum^{k_n}_{r=1}s_{j,r}$ is a sum of $jk_n$ i.i.d. random variables with distribution $\mbox{Exp}(1)$, that is, it has density $e^{-x}I(x\geq 0).$ Since the mean and the variance of $\mbox{Exp}(1)$ are both equal to $1$, we normalize the sum by
\beaa
W_j:=\frac{1}{\sqrt{jk_n}}\Big(\big(\sum^{k_n}_{r=1}s_{j,r}\big)-jk_n\Big).
\eeaa
By Theorem 1 on page 217 from Petrov (1975),
for any sequence of positive numbers $\delta_n$ such that $\delta_n=o((nk_n)^{1/6})$,
\begin{equation}\label{normapp}
P(W_j>x)=(1+o(1))(1-\Phi(x))~~~
\end{equation}
uniformly over $x\in [0,\delta_n]$ and $n/2\le j\le n$  as $n\to\infty$.
Now reorganize the index in (\ref{RiemannGeo}) to obtain
\bea
& & P\Big(T_n(j_n)\le k_n\psi(n)+\Big(\frac{k_n}{n}\Big)^{1/2}(a_n+b_ny)\Big) \nonumber\\
& = &   \prod^{j_n}_{i=1}P\Big(W_{n-i+1} \leq \big((n-i+1)k_n\big)^{1/2}(\psi(n)-\psi(n-i+1)) \nonumber\\
& & ~~~~~~~~~~~~~~~~~~~~~~~~~~~~~~~~~~~~~~~~~~~~~~~~ + \Big(\frac{n-i+1}{n}\Big)^{1/2}(a_n+b_ny)\Big) \nonumber\\
& =& \prod^{j_n}_{i=1}(1-a_{ni}),\lbl{baba}
\eea
where $a_{ni}=P\big(W_{n-i+1} > x_{n,i}\big)$ and
\bea\lbl{BowlL}
x_{n,i}=((n-i+1)k_n)^{1/2}(\psi(n)-\psi(n-i+1))+\Big(1-\frac{i-1}{n}\Big)^{1/2}
(a_n+b_ny).
\eea
Recalling (\ref{onjn}), we know $j_n=o(n)$. From the second expression in \eqref{psi}
we have
\[\psi(n)-\psi(n-i+1)=\frac{i-1}{n}\Big(1+O\Big(\frac{i}{n}\Big)\Big)
\]
uniformly over $1\le i\le j_n$ as $n\to\infty$. It follows that
\beaa
& & ((n-i+1)k_n)^{1/2}(\psi(n)-\psi(n-i+1))\\
&=&\Big(\frac{k_n}{n}\Big)^{1/2}(i-1)\cdot \Big(1-\frac{i-1}{n}\Big)^{1/2}\cdot \Big(1+O\Big(\frac{i}{n}\Big)\Big)\\
& = & \Big(\frac{k_n}{n}\Big)^{1/2}(i-1)\cdot \Big(1+O\Big(\frac{i}{n}\Big)\Big)\\
& = & \Big(\frac{k_n}{n}\Big)^{1/2}(i-1)\cdot \Big(1+O\big(n^{-3/8}\big)\Big)
\eeaa
uniformly over $1\le i\le j_n$ as $n\to\infty$. Since $a_n+b_ny=O((\log n)^{1/2})$,
\beaa
\Big(\Big(1-\frac{i-1}{n}\Big)^{1/2}-1\Big)(a_n+b_ny)
&=& (i-1)\cdot O\Big(\frac{(\log n)^{1/2}}{n}\Big)\\
&=& \Big(\frac{k_n}{n}\Big)^{1/2}(i-1)\cdot O\Big(\frac{(\log n)^{1/2}}{n^{1/2}}\Big)
\eeaa
uniformly over $1\le i\le j_n$.  Therefore, by combining the above two expansions we get
\begin{equation}\label{xn}
x_{n,i}=\Big(1+O\big(n^{-3/8}\big)\Big)
\Big(\frac{k_n}{n}\Big)^{1/2}(i-1)+a_n+b_ny
\end{equation}
uniformly over $1\le i\le j_n$. We emphasize the above is true when $i=j_n=1$, which can be seen directly from (\ref{BowlL}). This fact will be used later.

Finally, we prove \eqref{stepthree} by considering the three cases: $\alpha=0$, $\alpha\in (0, \infty)$ and $\alpha=\infty$.

{\it Case 1: $\alpha=0$}. Since
\beaa
& & a_n=a(n/k_n)\sim (\log (n/k_n))^{1/2}\ \ \mbox{and}\ \  b_n=b(n/k_n)=(\log (n/k_n))^{-1/2}\to 0,
\eeaa
we have
\[
\min_{1\le i\le j_n}x_{n,i}\to\infty\mbox{ and } \max_{1\le i\le j_n}x_{n,i}=O\Big(\frac{k_n^{1/2}j_n}{n^{1/2}}+(\log n)^{1/2}\Big)=O\big(n^{1/8}\big).
\]
It follows from \eqref{normapp} that
\[
a_{ni}=(1+o(1))(1-\Phi(x_{n,i}))
\]
uniformly over $1\le i\le j_n$. In Lemma~\ref{seed}, choose $x_n=n/k_n$, $j_n$ as in (\ref{onjn}) and $c_{nj}=(1+O(n^{-3/8}))
\big(\frac{k_n}{n}\big)^{1/2}$ as in (\ref{xn}) to obtain
\[
\sum^{j_n}_{i=1}a_{ni}=(1+o(1))\sum^{j_n}_{i=1}(1-\Phi(x_{n,i}))\to e^{-y}.
\]
Further, it is easily seen that $\max_{1\le i\le j_n}a_{ni}\to 0$. Applying (\ref{lem2}) to (\ref{baba}), we arrive at
\[
 P\Big(\frac{T_n(j_n)-k_n\psi(n)}{(\frac{k_n}{n})^{1/2}}\le  a_n+b_n y\Big)\to e^{-e^{-y}}=\Phi_0(y),\ \ y\in \mathbb{R},
\]
that is, we get \eqref{stepthree} for $\alpha=0$.

{\it Case 2}:  We see that  $j_n\sim \alpha^{-1/2}n^{1/8}$ from
\eqref{onjn}. By definition, $a_n=0$ and $b_n=1$.
Then it follows from  \eqref{xn} that
\[
x_{n,i}=(1+o(1))\alpha^{1/2}(i-1)+y
\]
holds uniformly over $1\le i\le j_n$ as $n\to\infty$. We claim that
\begin{equation}\label{candy}
a_{ni}=(1+o(1))(1-\Phi(x_{n,i}))
\end{equation}
uniformly over $1\le i\le j_n$.  In fact, review that \eqref{normapp} holds if $0<x=o(n^{1/3})$. Evidently, $\max_{1\leq i \leq j_n}|x_{n,i}|=O(n^{1/8}).$ But there is a possibility that $x_{n,i}<0$ for small values of $i$. Let $j_0>1$ be an integer such that  $\min_{j_0\le i\le j_n}x_{n,i}>0$. Then we
have from  \eqref{normapp} that \eqref{candy} holds uniformly over $j_0\le i\le j_n$. By using the standard central limit theorem, we know
\eqref{normapp} holds as well for each $i=1,\cdots, j_0-1$. Therefore, for each  $i\ge 1$,
\bea\lbl{heada}
\lim_{n\to\infty}a_{ni}= 1-\Phi\big(\alpha^{1/2}(i-1)+y\big)\ \ \ \mbox{and}\ \ \ \sum_{i\geq 1} \big(1-\Phi\big(\alpha^{1/2}(i-1)+y\big)\big) <\infty
\eea
by the fact $1-\Phi(x)\sim \frac{1}{\sqrt{2\pi}x}e^{-x^2/2}$ as $x\to+\infty$. We now apply Lemma~\ref{lem1} to show \eqref{stepthree}.
By defining $a_{ni}=0$ for all $i>j_n$, with \eqref{heada}, we only need to verify the following two conditions: $\sup_{n\ge n_0, 1\le i\le j_n}a_{ni}<1$ for some integer
$n_0$ and $\lim_{n\to\infty}\sum^{j_n}_{i=1}a_{ni}=\sum^{\infty}_{i=1}(1-\Phi(\alpha^{-1/2}(i-1)+y))$.
The first one follows from \eqref{candy} and the fact that $x_{n,i}\ge \frac{1}{2}\alpha^{1/2}(i-1)+y\ge y$ for
$1\le i\le j_n$ for all large $n$.
The second condition can be easily verified by the dominated convergence theorem since $a_{ni}\le 2(1-\Phi(\frac{1}{2}\alpha^{1/2}(i-1)+y))$ for all
$1\le i\le j_n$ as $n$ is sufficiently large and
$\sum^\infty_{i=1}2(1-\Phi(\frac{1}{2}\alpha^{1/2}(i-1)+y)<\infty$.

{\it Case 3: $\alpha=\infty$}. From \eqref{onjn}, $0\le (\frac{k_n}{n})(j_n-1)\le n^{1/8}$  and thus
$x_{n,i}=O(n^{1/8})$ by (\ref{xn}). In particular, we have $x_{n,1}=y$, and for all large $n$, $x_{n,i}>0$ if $2\le i\le j_n$ and $j_n\geq 2$.
Therefore, 
\beaa
a_{ni}=(1+o(1))(1-\Phi(x_{n,i}))
\eeaa
uniformly over $1\le i\le j_n$ from \eqref{normapp}. From (\ref{xn}), $a_{n1}\to 1-\Phi(y)$ as $n\to\infty$. Obviously, $x_{n,i}\geq \frac{i}{3}(\frac{k_n}{n})^{1/2}$ if $2\leq i\leq j_n$ and $j_n\geq 2$.   Thus, use the fact $1-\Phi(x)\sim \frac{1}{\sqrt{2\pi}\,x}e^{-x^2/2}$ as $x\to +\infty$ to see that, for large $n$,
\beaa
I(j_n\geq 2)\sum^{j_n}_{i=2}a_{ni}
&\leq & 2\sum^{\infty}_{i=2} \Big(1-\Phi\Big(\frac{i}{3}\Big(\frac{k_n}{n}\Big)^{1/2}\Big)\\
& \leq & \sum^{\infty}_{i=2} \exp\Big\{-\frac{k_n}{18n}i^2\Big\}\\
& \leq & \int_{-\infty}^{\infty}\exp\Big\{-\frac{k_n}{18n}x^2\Big\}\,dx = 3\sqrt{2\pi}\Big(\frac{n}{k_n}\Big)^{1/2}
\eeaa
since $\exp\big\{-\frac{k_n}{18n}i^2\big\} \leq \int_{i-1}^i\exp\big\{-\frac{k_n}{18n}x^2\big\}\,dx$ for all $i\geq 2.$ Thus,
$I(j_n\geq 2)\sum^{j_n}_{i=2}a_{ni}\to 0$. This and the fact $I(j_n\geq 2)\cdot\max_{2\leq i\leq j_n}a_{ni}\to 0$ imply that
$I(j_n\geq 2)\big(1-\prod^{j_n}_{i=2}(1-a_{ni})\big)\to 0$ as $n\to\infty$. So we have from (\ref{baba}) that
\beaa
& & P\Big(T_n(j_n)\le k_n\psi(n)+\Big(\frac{k_n}{n}\Big)^{1/2}(a_n+b_ny)\Big)\\
&=& \prod^{j_n}_{i=1}(1-a_{ni})\\
&=& (1-a_{n1})\cdot \Big[1+I(j_n\geq 2)\big(-1+\prod^{j_n}_{i=2}(1-a_{ni})\big)\Big]\to \Phi(y)=\Phi_{\infty}(y)
\eeaa
as $n\to\infty$. Reviewing the notation of $T_n(j_n)$ defined above (\ref{noxing}), we get \eqref{stepthree} for the case $\alpha=\infty$. The proof of the proposition is then completed.   \hfill $\blacksquare$\\


\noindent\textbf{Proof of Theorem \ref{Theorem_product}}. We use the same notation as in Proposition \ref{prop5}.
We first show the following:

\noindent (i) If $\lim_{n\to\infty}k_n/n=0$, particularly for $k_n\equiv k$, then
\bea\lbl{sleeped}
\frac{2(n/k_n)^{1/2}}{b_n}\Big(\frac{\max_{1\le j\le n}|z_j|}{n^{k_n/2}}-1\Big)-\frac{a_n}{b_n}
\ \ \mbox{converges weakly to cdf}\ \  \exp(-e^{-x}).
\eea

\noindent (ii) If $\lim_{n\to\infty}k_n/n=\alpha\in (0,\infty)$, then
\bea\lbl{nono}
\frac{\max_{1\le j\le n}|z_j|}{n^{k_n/2}}\ \ \mbox{converges weakly to cdf}\ \  \Phi_{\alpha}\Big(\frac{1}{2}\alpha^{1/2}+2\alpha^{-1/2}\log x\Big),~~~~x>0.
\eea
To do so,  for $\alpha\in [0,\infty)$, define
\[
V_n=\frac{\max_{1\le j\le n}\log|z_j|-k_n\psi(n)/2}{(k_n/n)^{1/2}b_n/2}-\frac{a_n}{b_n}.
\]
Then $V_n$ converges in distribution
to $\Theta_\alpha$ by Proposition \ref{prop5}, where $\Theta_\alpha$ is a random variable with cdf $\Phi_\alpha(y)$. Trivially,
\begin{eqnarray}\label{zzz}
\max_{1\le j\le n}|z_j|&=&\exp\Big\{\frac12k_n\psi(n)+\frac12\Big(\frac{k_n}{n}\Big)^{1/2}(a_n+b_nV_n)\Big\}\nonumber\\
&=&\exp\Big\{\frac12\Big(k_n\psi(n)+\Big(\frac{k_n}{n}\Big)^{1/2}a_n\Big)\Big\}\cdot\exp\Big\{\frac12\Big(\frac{k_n}{n}\Big)^{1/2}b_nV_n\Big\}.
\end{eqnarray}

If $\alpha=0$, then $\frac{k_n}{n}\to 0$,  $a_n=a(\frac{n}{k_n})\sim (\log\frac{n}{k_n})^{1/2}\to\infty$,
$b_n=b(\frac{n}{k_n})=(\log\frac{n}{k_n})^{-1/2}\to 0$,  and
$(\frac{k_n}{n})^{1/2}a_n\sim (\frac{k_n}{n})^{1/2}b_n^{-1}\to 0$ as $n\to\infty$. Using \eqref{psi} and expanding
\eqref{zzz} we get
\begin{eqnarray*}
\max_{1\le j\le n}|z_j|
&=&\exp\Big\{\frac12k_n\log n+O\Big(\frac{k_n}{n}\Big)+\frac12\Big(\frac{k_n}{n}\Big)^{1/2}a_n\Big\}\Big(1+\frac12\Big(\frac{k_n}{n}\Big)^{1/2}b_nV_n
+O_P\Big(\frac{k_nb_n^2}{n}\Big)\Big)\\
&=&n^{k_n/2}\Big(1+\frac12\Big(\frac{k_n}{n}\Big)^{1/2}a_n+O\Big(\frac{k_n}{n}\Big)\Big)\Big(1+\frac12\Big(\frac{k_n}{n}\Big)^{1/2}b_nV_n
+O_P\Big(\frac{k_n}{n}\Big)\Big)\\
&=&n^{k_n/2}\Big(1+\frac12\Big(\frac{k_n}{n}\Big)^{1/2}a_n+\frac12\Big(\frac{k_n}{n}\Big)^{1/2}b_nV_n
+O_P\Big(\frac{k_na_n^2}{n}\Big)\Big),
\end{eqnarray*}
 which yields that
 \[
 \frac{2(n/k_n)^{1/2}}{b_n}\Big(\frac{\max_{1\le j\le n}|z_j|}{n^{k_n/2}}-1\Big)-\frac{a_n}{b_n}
 =V_n+O_P\Big(\Big(\frac{k_n}{n}\Big)^{1/2}\Big(\log\frac{n}{k_n}\Big)^{3/2}\Big)
\]
converges in distribution to $\Lambda$ by the Slutsky lemma. We obtain (\ref{sleeped}).

Now assume $\alpha\in (0,\infty)$. In this case, $a_n=0$ and $b_n=1.$ Then from \eqref{zzz},
\[
\max_{1\le j\le n}|z_j|
=\exp\Big\{\frac12k_n\psi(n)\Big\}\exp\Big\{\frac12\Big(\frac{k_n}{n}\Big)^{1/2}V_n\Big\}.
\]
Using expansion $\psi(n)=\log n-\frac{1}{2n}+O(\frac{1}{n^2})$ from \eqref{psi}  we have
\[
\frac{\max_{1\le j\le n}|z_j|}{n^{k_n/2}}=\exp\Big\{-\frac14\alpha+o(1)\Big\}\cdot\exp\Big\{\Big(\frac12\alpha^{1/2}+o(1)\Big)V_n\Big\},
\]
which converges weakly to the distribution of  $e^{-\alpha/4}\exp\big(\frac12\alpha^{1/2}\Theta_{\alpha}\big)$, given by $\Phi_{\alpha}(\frac12\alpha^{1/2}+ 2\alpha^{-1/2}\log y)$, $y>0$.  We get (\ref{nono}).

From (\ref{sleeped}) it is easy to see
\beaa
\frac{2(n/k_n)^{1/2}}{b_n}= \Big(\frac{n}{k_n}\log \frac{n}{k_n}\Big)^{1/2}=\alpha_n\ \ \ \mbox{and}\ \ \ \frac{a_n}{b_n}=\log \frac{n}{k_n}-\log \log \frac{n}{k_n}-\frac{1}{2}\log (2\pi)=\beta_n.
\eeaa
Thus we obtain (a) of Theorem \ref{Theorem_product}. The part (b) follows from (\ref{nono}) and the part (c) is yielded from Proposition \ref{prop5} with $\Phi_\infty(x)=\Phi(x)$. This completes the proof of the theorem. \hfill$\blacksquare$\\

\subsection{The Verifications of (\ref{smartyy}) and (\ref{wisconsin})}\lbl{Flyingcolor}

\noindent\textbf{Verification of (\ref{smartyy})}. First, by the Taylor expansion,
\beaa
H_k(y)=e^{-y}\sum_{j=0}^{k-1}\frac{y^j}{j!} &=& 1-e^{-y}\sum^{\infty}_{j=k}\frac{y^j}{j!}\\
&=& 1-\frac{y^ke^{-y}}{k!}\Big(1+ y\sum^{\infty}_{j=k+1}\frac{y^{j-k-1}}{(k+1)\cdots j}\Big)
\eeaa
for all $y\in \mathbb{R}$ and $k\geq 1.$ Notice that the absolute value of the above sum is bounded by  $\sum^{\infty}_{j=k+1}\frac{1}{(k+1)\cdots j}\leq 1+\sum_{j=1}^{\infty}\frac{1}{j^2}<\infty$ uniformly for all $k\geq 1$ and $|y|\leq 1$. This says
\beaa
H_k(y)=1-\frac{y^k}{k!}e^{-y}(1+O(y))=1-\frac{y^k}{k!}(1+O(y))
\eeaa
as $y\to 0$ uniformly for all $k\geq 1.$ Hence
\beaa
\log \prod^\infty_{k=1}H_k(y) & = &  \sum^\infty_{k=1}\log H_k(y)\\
&=& -(1+O(y))\sum^\infty_{k=1}\frac{y^k}{k!} =-y(1+O(y))
\eeaa
since $\sum^\infty_{k=1}\frac{y^k}{k!}=e^y-1 \sim y$ as $y\to 0.$ Therefore,
\beaa
1-\prod^\infty_{k=1}H_k(y) =1- e^{-y(1+O(y))} \sim y
\eeaa
as $y\to 0.$ Taking $y=x^{-2}$ and letting $x\to \infty$, we get (\ref{smartyy}). \hfill$\blacksquare$\\

\noindent\textbf{Verification of (\ref{wisconsin})}. Given parameter $\beta>0,$ set
\[
F_\beta(x)=\prod^\infty_{j=0}\Phi\Big(x+\beta j\Big)
\]
for $x\in \mathbb{R}.$ From integration by parts, we know $1-\Phi(x)=\frac{1}{\sqrt{2\pi}\, x}e^{-x^2/2}(1+O(x^{-2}))$ as $x\to +\infty.$ Use $\log (1-t)=-t(1+O(t))$ as $t\to 0$ to have
\beaa
\log \Phi(x) & = & -\frac{1}{\sqrt{2\pi}\, x}e^{-x^2/2}(1+O(x^{-2}))\Big(1+ O\Big(\frac{1}{x}e^{-x^2/2}\Big)\Big)\\
& = & -\frac{1}{\sqrt{2\pi}\, x}e^{-x^2/2}(1+a(x))
\eeaa
as $x\to+\infty$, where $a(x)$ is defined over $[1,\infty)$ and $|a(x)|\leq Cx^{-2}$ for all $x\geq 1$ and $C$ is a constant not depending on $x$. Thus,
\bea
\log F_\beta(x) & = & \sum^\infty_{j=0}\log \Phi\Big(x+\beta j\Big) \nonumber\\
& = & -\frac{1}{\sqrt{2\pi}}\sum^\infty_{j=0} \frac{1}{x+\beta j}e^{-(x+\beta j)^2/2}(1+a(x+\beta j))\nonumber\\
& = & -(1+o(1))\frac{1}{\sqrt{2\pi}}\sum^\infty_{j=0} \frac{1}{x+\beta j}e^{-(x+\beta j)^2/2} \lbl{sound_la}
\eea
as $x\to +\infty$. Observe
\beaa
\frac{1}{x+\beta (j+1)}e^{-(x+\beta (j+1))^2/2} \leq \int_{x+\beta j}^{x + \beta (j+1)}\frac{1}{t}e^{-t^2/2}\,dt \leq \frac{1}{x+\beta j}e^{-(x+\beta j)^2/2}
\eeaa
for all $x>0$ and $j\geq 0$. Sum the above over all $j\geq 1$ to obtain
\beaa
\int_{x+\beta}^{\infty}\frac{1}{t}e^{-t^2/2}\,dt\leq \sum^\infty_{j=1} \frac{1}{x+\beta j}e^{-(x+\beta j)^2/2} \leq \frac{1}{x+\beta}e^{-(x+\beta)^2/2} + \int_{x+\beta}^{\infty}\frac{1}{t}e^{-t^2/2}\,dt
\eeaa
for all $x>0.$ Write $\int_{x+\beta}^{\infty}\frac{1}{t}e^{-t^2/2}\,dt=-\int_{x+\beta}^{\infty}\frac{1}{t^2}(e^{-t^2/2})'\,dt$. From the integration by parts, $\int_{x+\beta}^{\infty}\frac{1}{t}e^{-t^2/2}\,dt \sim \frac{1}{(x+\beta)^2}e^{-(x+\beta)^2/2}$ as $x\to+\infty$. Since $\beta>0$, we have $\frac{1}{(x+\beta)^2}e^{-(x+\beta)^2/2}=o\big(\frac{1}{x}e^{-x^2/2}\big)$ and $\frac{1}{x+\beta}e^{-(x+\beta)^2/2}=o\big(\frac{1}{x}e^{-x^2/2}\big)$ as $x\to +\infty.$ It follows from (\ref{sound_la}) that
\beaa
\log F_\beta(x)
\sim -\frac{1}{\sqrt{2\pi}}\frac{1}{x}e^{-x^2/2}
\eeaa
as $x\to +\infty$. In other words, the first term in the sum  appeared in (\ref{sound_la}) dominates the sum.   Thus,
\beaa
1-F_\beta(x)=1-e^{\log F_\beta(x)}\sim -\log F_\beta(x) \sim \frac{1}{\sqrt{2\pi}}\frac{1}{x}e^{-x^2/2}
\eeaa
as $x\to+\infty.$ Observe that the above approximation is free of the choice of $\beta$. Since $F_{\sqrt{\alpha}}(x)=\Phi_{\alpha}(x)$ for $x>0.$
Replacing ``$x$" by ``$\frac{1}{2}\alpha^{1/2}+2\alpha^{-1/2}\log x$", we arrive at
\beaa
& & 1-\Phi_{\alpha}\big(\frac{1}{2}\alpha^{1/2}+2\alpha^{-1/2}\log x\big)\\
 & = & 1-F_{\sqrt{\alpha}}\big(\frac{1}{2}\alpha^{1/2}+2\alpha^{-1/2}\log x\big)\\
& \sim & \frac{1}{\sqrt{2\pi}}\big(\frac{1}{2}\alpha^{1/2}+2\alpha^{-1/2}\log x\big)^{-1}\exp\Big\{-\big(\frac{1}{2}\alpha^{1/2}+2\alpha^{-1/2}\log x\big)^2/2\Big\}\\
& \sim & \frac{\sqrt{\alpha} e^{-\alpha/8}}{2\sqrt{2\pi}}\frac{1}{x\log x}e^{-2(\log x)^2/\alpha}
\eeaa
as $x\to +\infty$. At last
\bea
P\big(e^{N(0,1)}\geq x\big)=P(N(0,1)\geq \log x) \sim \frac{1}{\sqrt{2\pi}\,\log x}e^{-(\log x)^2/2}
\eea
as $x\to +\infty$. This  verifies (\ref{wisconsin}) and the statement below.   \hfill$\blacksquare$\\

\noindent\textbf{Acknowledgements}. We thank Drs. Ming Gao, Wenqing Hu, Jing Wang, Ke Wang and Gongjun Xu for helping us check the proofs.

\baselineskip 12pt
\def\ref{\par\noindent\hangindent 25pt}

\end{document}